\documentclass[12pt]{smfart}
\usepackage{amssymb}
\usepackage{amscd}
\usepackage{amsmath}
\usepackage[french]{babel}
\usepackage{t1enc}
\usepackage{mathrsfs}
\usepackage[all]{xypic}

\def\cartesien{\ar@{}[rd]|{\square}}

\newcommand{\PGL}{\ensuremath{\operatorname{PGL}}}

\newcommand{\CH}{\ensuremath{\operatorname{CH}}}
\newcommand{\chow}{\ensuremath{\operatorname{CH}}}

\newcommand{\chowtors}[2]{\ensuremath{\operatorname{CH}^{#1}(#2)_{\operatorname{tors}}}}

\newcommand{\SB}[2]{\ensuremath{\operatorname{SB}_{#1}(#2)}}

\newcommand{\Flag}[2]{\ensuremath{\operatorname{SB}_{#1}(#2)}}

\newcommand{\Drapeau}[2]{\ensuremath{\operatorname{Drap}_{#1}(#2)}}

\newcommand{\Spec}[1]{\ensuremath{\operatorname{Spec}(#1)}}

\renewcommand{\Im}{\ensuremath{\operatorname{Im}}}
\newcommand{\ind}{\ensuremath{\operatorname{ind}}}
\newcommand{\End}{\ensuremath{\operatorname{End}}}

\newcommand{\Hom}{\ensuremath{\operatorname{Hom}}}
\newcommand{\brauer}{\ensuremath{\operatorname{Br}}}
\newcommand{\Pic}{\ensuremath{\operatorname{Pic}}}

\newcommand{\enscol}[2]{\ensuremath{\left\lbrace\begin{array}{c}{#1}\\{#2}\end{array}\right\rbrace}}

\newcommand{\enscolquatre}[4]{\ensuremath{\left\lbrace\begin{array}{c}{#1}\\{#2}\\{#3}\\{#4}\end{array}\right\rbrace}}

\theoremstyle{plain}
\newtheorem{thm}{Th\'eor\`eme}[section]
\newtheorem{lem}[thm]{Lemme}
\newtheorem{prop}[thm]{Proposition}
\newtheorem{cor}[thm]{Corollaire}

\theoremstyle{definition}
\newtheorem{defi}[thm]{D\'efinition}

\newtheorem{nota}[thm]{Notation}

\theoremstyle{remark}
\newtheorem{rem}[thm]{Remarque}

\renewcommand{\proofname}{proof}
\def\noqed{\let\qed\relax}

\makeatletter
\def\markboth#1#2{%
  \begingroup
    \@temptokena{{#1}{#2}}\xdef\@themark{\the\@temptokena}%
    \mark{\the\@temptokena}%
  \endgroup
  \if@nobreak\ifvmode\nobreak\fi\fi}
\makeatother

\author{Franck Doray}
\address{Mathematisch Instituut\\
Universiteit Leiden\\
Postbus 9512\\
2300 RA Leiden}
\email{doray@math.leidenuniv.nl}

\title{Vari\'et\'es homog\`enes sous $\PGL_n$}

\begin{document}

\begin{altabstract}
Let $A$ be an Azumaya algebra over a field. If $G$ is the group of automorphisms of $A$ and $X$ denotes a projective homogeneous variety under $G$, we construct in a very explicit way and under suitable hypotheses a bundle $\mathcal{V}$ on $S$, where $S$ is a (generalized) Severi-Brauer variety associated to $A$, and a canonical isomorphism between $X$ and a flag bundle on $\mathcal{V}$. This allows to explicitely compute Chow groups of $X$ in terms of the Chow groups of $S$. 
\end{altabstract}

\begin{abstract}
Soit $A$ une alg\`ebre d'Azumaya sur un corps. Notons $G$ le groupe de ses automorphismes. 
Si $X$ d\'esigne une vari\'et\'e homog\`ene projective sous $G$, nous construisons explicitement, sous certaines hypoth\`eses,  un fibr\'e $\mathcal{V}$ sur $S$, o\`u $S$ est une vari\'et\'e de Severi-Brauer (g\'en\'eralis\'ee) associ\'ee \`a $A$, et un isomorphisme canonique entre $X$ et une fibration en drapeaux de $\mathcal{V}$. Ceci permet de calculer de mani\`ere explicite les groupes de Chow de X en fonction des groupes de Chow de $S$. 
\end{abstract}
\maketitle

\section{Introduction}
Si  une $k$-vari\'et\'e projective $X$ est suppos\'ee cellulaire, alors l'anneau de Chow gradu\'e $\bigoplus_i \chow^i (X)$ est le $\mathbf{Z}$-module libre sur les classes des adh\'erences des cellules de $X$ (\cite{kahn_motivic}). En particulier, les vari\'et\'es de drapeaux d\'eploy\'ees,
c'est-\`a-dire les vari\'et\'es projectives homog\`enes sous un groupe alg\'ebrique semi-simple d\'eploy\'e ont des groupes
de Chow libres de type fini sur $\mathbf{Z}$. Il est alors l\'egitime de vouloir \'etudier les groupes de Chow de vari\'et\'es projectives homog\`enes non d\'eploy\'ees. Consid\'erons une telle $k$-vari\'et\'e projective homog\`ene $X$ et supposons-la d\'eploy\'ee sur une extension $L/k$. Alors la fl\`eche naturelle $\chow^i (X) \to \chow^i (X_L)$ permet de voir que tout \'el\'ement de torsion de $\chow^i (X)$ s'envoie sur 0, et en fait le noyau est exactement constitu\'e d'\'el\'ements de torsion (la formule de projection permet en effet de voir que la compos\'ee $\chow^i (X) \to \chow^i (X_L) \to \chow^i (X)$ co\"incide avec la multiplication par le degr\'e de l'extension $L/k$).

Les groupes de Chow en codimension 1 sont connus (voir par exemple \cite[lemme 5.1]{peyre_galois3}). N. Karpenko a calcul\'e les groupes de Chow en codimension 2 des vari\'et\'es de Severi-Brauer, cas particulier de vari\'et\'es homog\`enes projectives sous $\PGL(A)$ pour une alg\`ebre simple centrale $A$. Le calcul  des motifs des  vari\'et\'es homog\`enes projectives sous un groupe $G$, forme interne d'un groupe de Chevalley se ram\`ene dans certains cas au calcul de motifs de  vari\'et\'es homog\`enes projectives \og plus simples \fg (i.e.: $G/P$ o\`u $P$ est un parabolique maximal) comme l'ont montr\'e Calm\`es, Petrov, Semenov et Zainouline dans \cite{calmes_motives} . Le pr\'esent article a pour but de calculer de mani\`ere tout-\`a-fait explicite les groupes de Chow de  vari\'et\'es homog\`enes projectives sous $\PGL_n$ en fonction des groupes de Chow de vari\'et\'es de Severi-Brauer ou Severi-Brauer g\'en\'eralis\'ees (\cite[d\'efinition 1.16]{invol}).

\section{Quelques notations et d\'efinitions}\label{nota_defi}

Introduisons ici quelques notations et d\'efinitions que nous utiliserons par la suite.

Nous utiliserons la th\'eorie des alg\`ebres simples centrales que nous appellerons aussi alg\`ebres d'Azumaya~:
\index{Alg\`ebre simple centrale}
\index{Alg\`ebre d'Azumaya}
\begin{defi}
Notons $k_s$ une cl\^oture s\'eparable de $k$.
Soit $A$ une $k$-alg\`ebre. On dira que $A$ est une alg\`ebre d'Azumaya s'il existe un entier $r\geqslant 1$ tel que
$$
A\otimes_k k_s = \mathcal{M}_r (k_s).
$$
\end{defi}

Une telle alg\`ebre est de dimension $r^2$, $r$ est appel\'e le degr\'e de $A$.

Si $X$ est un sch\'ema, rappelons \cite{groth_dix_1} qu'une alg\`ebre $\mathcal{A}$ sur $X$ est une alg\`ebre d'Azumaya si, localement pour la topologie \'etale, $\mathcal{A}$ est isomorphe \`a $\mathcal{M}_r(\mathcal{O}_X)$ pour un certain $r$. Si $X=\Spec{k}$ est le spectre d'un corps, une alg\`ebre d'Azumaya $A$ sur $X$ n'est rien d'autre qu'une alg\`ebre simple centrale sur $k$.

\begin{defi}\label{def-flag}
Soient $A$ une alg\`ebre d'Azumaya sur $k$ de degr\'e $n$ et $1\leqslant i_1<\ldots<i_r \leqslant n$ des entiers. On notera par $\Flag{i_1,\ldots,i_r}{A}$ le $k$-sch\'ema repr\'esentant le foncteur qui \`a une $k$-alg\`ebre $R$ associe l'ensemble des $r$-uplets $( I_1,\ldots,I_r )$ d'id\'eaux \`a gauche de $A\otimes_k R$, tels que $A\otimes_k R /I_j$ soit localement libre de rang $n^2-ni_j$ et $I_1\subset\ldots \subset I_r$.
\end{defi}

\begin{rem}
Tout groupe alg\'ebrique $G$, forme int\'erieure de type $A_{n-1}$ s'\'ecrit $\PGL(A)$ pour une alg\`ebre d'Azumaya $A$ de degr\'e $n$, et toute vari\'et\'e projective homog\`ene $X$ sur $G$ s'\'ecrira $\Flag{i_1,\ldots,i_r}{A}$ pour un entier $r$ et un $r$-uplet $(i_1,\ldots,i_r)$  \cite[5.2]{merkurjev_index_1}.
\end{rem}

\begin{rem}
En particulier $\Flag{i}{A}$ est une vari\'et\'e de Severi-Brauer g\'en\'eralis\'ee 
\cite[d\'efinition 1.16]{invol}.
Par exemple, $\Flag{1}{A}$ co\"incide avec $\SB{}{A}$, vari\'et\'e de Severi-Brauer \og habituelle \fg , d\'efinie pour la premi\`ere fois par Ch\^atelet \cite{chat}.      
\end{rem}

\begin{rem}
Si $L/k$ est une extension de corps  d\'eployant $A$ (un tel $L$ existe d'apr\`es le th\'eor\`eme de Wedderburn, \cite[th\'eor\`eme 1.1]{invol}), alors $\Flag{i}{A}\otimes_{\Spec{k}} \Spec {L}$ est simplement une vari\'et\'e de Grassmann, d'apr\`es \cite[th\'eor\`eme 1.18]{invol}.
\end{rem}

\begin{nota}\label{nota_drapeaux}
Soit $V$ un $k$-espace vectoriel de dimension $n$. Soient $1\leqslant i_1<\ldots<i_r \leqslant n$ des entiers. La vari\'et\'e de drapeaux sur $V$ associ\'ee \`a la suite $(1\leqslant i_1<\ldots<i_r \leqslant n)$ est not\'ee
$\Drapeau{i_1,\ldots,i_r}{V}$.  C'est le $k$-sch\'ema repr\'esentant le foncteur qui \`a une $k$-alg\`ebre $R$ associe l'ensemble des $r$-uplets $(V_1,\ldots,V_r)$ de $R$-modules de $V\otimes_k R$ tels que $V\otimes_k R /V_i$ soit localement libre de rang $n-i_j$ et $V_1 \subset \ldots \subset V_r$. Par exemple, ses $L$-points pour une extension de corps $L/k$ sont $\Drapeau{i_1,\ldots,i_r}{V}(L)=\lbrace V_{1} \subset \ldots \subset V_{r}\subset V\otimes_k L, \mbox{ pour tout $j$, $V_j$, $L$-espace vectoriel de dimension $i_j$}\rbrace$.
\end{nota}

\begin{defi}\label{fibration-drapeaux}
Soient $\mathcal{V}$  un fibr\'e de rang $n$ sur un $k$-sch\'ema $X$, et $1\leqslant i_1<\ldots<i_r \leqslant n$ une suite d'entiers. On notera $\Drapeau{i_1,\ldots,i_r}{\mathcal{V}}$ le fibr\'e en vari\'et\'es de drapeaux associ\'e \`a $1\leqslant i_1<\ldots<i_r \leqslant n$ d\'efini comme suit:
si $U$ est un ouvert de $X$ trivialisant $\mathcal{V}$, on a le diagramme commutatif suivant dont le carr\'e est cart\'esien :
$$
\xymatrix{
   \mathbf{A}^n_k \times U 
   \ar[dr]_-{pr_2}
   &
   \mathcal{V}_{|U} 
   \ar[l]^-{\sim}_-\phi
   \ar[r]
   \ar[d]
   \cartesien
   &
   \mathcal{V}
   \ar[d]
   \\
   &
   U
   \ar[r]
   &
   X}                                      
$$

\noindent
On d\'efinit un fibr\'e en vari\'et\'es de drapeaux $\mathcal{E}_U$ sur $U$ comme \'etant \'egal \`a $\phi^{-1}(\Drapeau{i_1,\ldots,i_r}{k^n}\times U)$. Pour tout ouvert $U$ trivialisant $\mathcal{V}$, on construit ainsi un fibr\'e en vari\'et\'es de drapeaux $\mathcal{E}_U$. Il est \'evident qu'ils se recollent pour donner naissance \`a un fibr\'e en vari\'et\'es de drapeaux $\mathcal{E}$ sur $X$ tel que si $U$ est un ouvert trivialisant $\mathcal{V}$, on a : $\mathcal{E}_{|U}=\mathcal{E}_U$. Finalement, on pose : $\Drapeau{i_1,\ldots,i_r}{\mathcal{V}}=\mathcal{E}$.
\end{defi}

\begin{nota}
Si $\mathcal{B}$ est une alg\`ebre d'Azumaya de degr\'e $n$, sur un sch\'ema $X$, et si $1\leqslant i_1<\ldots <i_r \leqslant n$ sont des entiers, alors on peut d\'efinir une fibration en vari\'et\'es de drapeaux tordues $\Flag{i_1,\ldots,i_r}{\mathcal{B}} \to X$ qui revient \`a \'etendre la d\'efinition \ref{def-flag} \og en famille \fg ~(construction analogue \`a la d\'efinition \ref{fibration-drapeaux}). Par exemple, si $x\in X$ est un point, alors
$$\Flag{i_1,\ldots,i_r}{\mathcal{B}}_x=\Flag{i_1,\ldots,i_r}{\mathcal{B}_x}.$$
\noindent
Nous  noterons  $\SB{}{\mathcal{B}}=\SB{1}{\mathcal{B}}$.
\end{nota}

\section{R\'esultats}

Nous utiliserons les notations introduites en \ref{nota_defi}. \'Enon\c{c}ons le premier th\'eor\`eme:
\begin{thm}\label{thm_effectif}
Soit $A$ une alg\`ebre d'Azumaya sur $k$ de degr\'e $n$. Soient $r>1$ et $1 \leqslant i_1<\ldots<i_r \leqslant n$ des entiers. Supposons $i_1=1$, alors il existe un fibr\'e de rang $n-1$, $\mathcal{V}$ sur $\SB{i_1}{A}=\SB{}{A}$, tel que l'on ait un isomorphisme :
$$
\xymatrix{\Flag{i_1,i_2,\ldots,i_r}{A}\ar[rd]\ar[rr]^\sim&&\Drapeau{n-i_r,\ldots,n-i_2}{\mathcal{V}}\ar[ld]\\
&\SB{}{A}&}
$$
\end{thm}

Si $i_1$ est quelconque, on a le r\'esultat g\'en\'eral suivant :

\begin{thm}\label{thm_general}
Soit $A$ une alg\`ebre d'Azumaya sur $k$ de degr\'e $n$. Soient $r>1$ et $1 \leqslant i_1<i_2<\ldots<i_r \leqslant n$ des entiers. Fixons un $s\in \lbrace 1,\ldots, r \rbrace$. Il existe deux alg\`ebres d'Azumaya $\mathcal{B_{+}}$ et $\mathcal{B_{-}}$ sur $\SB{i_s}{A}$ de degr\'es respectifs $n-i_s$ et $i_s$ telles que l'on ait un isomorphisme :

$$
\xymatrix@-1pc{
  \Flag{i_1,\ldots,i_r}{A}
  \ar@/_/[rd]
  \ar[rr]^-\sim
  &&
  \Flag{i_1,\ldots,i_{s-1}}{\mathcal{B_{-}}}
  \times_{\SB{i_s}{A}}
  \Flag{i_{s+1}-i_s,\ldots,i_r-i_s}{\mathcal{B_{+}}}
  \ar@/_/[dl]
  \\
  &
  \SB{i_s}{A}
  &}
$$
\end{thm}

Il admet pour corollaires : 
\begin{cor}\label{cor_thm_general}
Conservons les notations du th\'eor\`eme \ref{thm_general}. Si l'indice de $A$, $\ind A$, (cf. notation \ref{nota_indice} ) est premier avec un des $i_j$, disons $i_s$, alors la projection naturelle~: 
$$\Flag{i_1,\ldots,i_r}{A} \to \SB{i_s}{A}$$ est une fibration en produits de vari\'et\'es de drapeaux non tordues.
\end{cor}

\begin{cor}\label{cor_ind_premier}
Conservons les notations du th\'eor\`eme \ref{thm_general}. Si $\ind A$ est une puissance d'un nombre premier et si $(\ind A,i_1,\ldots,i_r)=1$ alors il existe un $s\in \lbrace 1\ldots r \rbrace$ tel que la projection naturelle~: $$\Flag{i_1,\ldots,i_r}{A} \to \SB{i_s}{A}$$ soit une fibration en produits de vari\'et\'es de drapeaux non tordues.
\end{cor}

\section{Preuve des r\'esultats}

\Subsection{Espaces vectoriels sur un corps gauche - Annulateurs}\label{section_espacevectoriel}

Rappelons un r\'esultat classique de Wedderburn pour fixer les notations : 
\begin{thm}[Wedderburn]\label{wedd}
Si $A$ est une alg\`ebre d'Azumaya sur $k$ de degr\'e $n$, alors il existe un corps gauche $D$ sur $k$ et un $D$-espace vectoriel \`a droite $E$, tels que 
$$A=\End_D(E).$$
L'anneau des $D$-endomorphismes de $E$ op\`ere \`a gauche sur $E$.
\end{thm}
\begin{proof}
\cite[th\'eor\`eme 1.1]{invol}.
\end{proof}
Nous utiliserons librement ce r\'esultat dans la suite avec les m\^emes notations.
\begin{nota}\label{nota_indice}
La dimension de $D$ est n\'ecessairement un carr\'e, sa racine carr\'e sera l'indice de $A$, not\'e $\ind A$. Si $r$ est la $D$-dimension de $E$, alors le degr\'e de $A$ est \'egal \`a $rd$, o\`u $d=\ind A$.
\end{nota}

\begin{defi}
Soit $A$ une alg\`ebre d'Azumaya sur $k$, de degr\'e $n$.
Soit $I \subset A$, alors on d\'efinit son \textit{annulateur} \`a droite comme suit:
$I^{\circ}=\{a \in A,  I a =0\}$. De m\^eme, on d\'efinit son annulateur \`a gauche $^\circ \! I = \{ a \in A, a I = 0 \}$.
\end{defi}

\begin{rem}\label{rem_annulateurs}
Si $I$ est un sous-ensemble quelconque de $A$, alors $I^\circ$ est un id\'eal \`a droite de $A$ et $^\circ \! I$ est un id\'eal \`a gauche de $A$.
\end{rem}

La proposition classique suivante est d'une importance cruciale pour mieux comprendre les objets utilis\'es :

\begin{prop}\label{id-vect}
Soient $D$ un corps gauche de centre $k$, $E$ un $D$-espace vectoriel \`a droite de dimension finie et posons $A=\End_D(E)$, alg\`ebre d'Azumaya sur $k$.
L'application $V \to \Hom_D(E/V,E)$ d\'efinit une bijection entre, d'une part, les sous-$D$-espaces vectoriels de $E$ de $D$-dimension $l$ et, d'autre part, les id\'eaux \`a gauche de $A$ de $k$-dimension 
$\deg A(\deg A - l \ind A)$.
De m\^eme, $V \to \Hom_D(E,V)$ d\'efinit une bijection entre d'une part les sous-$D$-espaces vectoriels de $E$ de $D$-dimension $l$ et d'autre part les id\'eaux \`a droite de $A$ de $k$-dimension 
$l \deg A \ind A$.
\end{prop}
\begin{proof}
\cite[Proposition 1.12]{invol}. 
\end{proof}

\begin{rem}
Avec les notations pr\'ec\'edentes, si $d=\ind A=\deg D$, $n=\deg A$, et si $r$ est la $D$-dimension de $E$ ($n=dr$), alors le sous-$D$-espace vectoriel correspondant \`a un id\'eal \`a gauche $I$ de $A$ de $k$-dimension $ni$ (pour un certain $i$) est de $D$-dimension $\frac{n-i}{d}$.
\end{rem}

\begin{cor}\label{perp}
Conservons les notations de la proposition.
Si $I$ est un id\'eal \`a gauche de $A$ s'\'ecrivant $\Hom_D(E/V,E)$ pour un sous-$D$-espace vectoriel $V$ de $E$, alors $I^\circ$ s'identifie canoniquement \`a l'id\'eal $\Hom_D(E,V)$. De m\^eme si $I$ est un id\'eal \`a droite de $A$ s'\'ecrivant $\Hom_D(E,V)$ pour un sous-$D$-espace vectoriel $V$ de $E$, alors $^\circ \! I$ s'identifie canoniquement \`a l'id\'eal alors \`a gauche $\Hom_D(E/V,E)$.
\end{cor}
\begin{proof}
\cite[Preuve de la proposition 1.14.]{invol}
\end{proof}
On en d\'eduit :
\begin{lem}\label{lem_double_perp}
 Si $I$ est un id\'eal \`a gauche de $A$, alors $^\circ \!(I^\circ)=I$. De m\^eme si $I$ est un id\'eal \`a droite de $A$, alors $(^\circ \!I)^\circ =I$.
\end{lem}
\begin{proof}
 \cite[Proposition 1.14.]{invol}
\end{proof}

\begin{lem}\label{dimension}
Soit $A$ une alg\`ebre d'Azumaya sur $k$, de degr\'e $n$.\\
(i) Si $I \subset A$ est un id\'eal \`a gauche de $k$-dimension $ni$ alors $I^\circ$ est un id\'eal \`a droite de $k$-dimension $n(n-i)$.\\
(ii) Si $I$ et $J$ sont des id\'eaux \`a gauche de $A$ de $k$-dimensions respectives $ni$ et $nj$ alors $I^\circ J=I^\circ \cap J$ est un sous-$k$-espace vectoriel de $A$ de dimension $(n-i)j$.
\end{lem}
\begin{proof}
On peut \'ecrire : $A=\End_D(E)$.\\
(i) La remarque \ref{rem_annulateurs} permet de voir que $I^\circ$ est un id\'eal \`a droite; la proposition \ref{id-vect} nous donne sa dimension.\\
(ii) Si $I=\Hom_D(E/V,E)$ et $J=\Hom_D(E/W,E)$ o\`u $V$ et  $W$ sont des sous-$D$-espaces vectoriels de $E$, alors $I^\circ = \Hom_D(E,V)$, et un simple calcul permet de voir que $I^\circ \cap J =I^\circ J=\Hom_D(E/W,V)$. Par suite (proposition \ref{id-vect}) $\dim_k I^\circ J = \dim_k I^\circ \cap J=(n-i)j$.
\end{proof}

Le lemme suivant est \`a rapprocher des propositions 1.15 et 1.20 de \cite{invol}.
\begin{lem}\label{quot1}
Conservons les notations de la proposition \ref{id-vect}. Si $I$ est un id\'eal \`a gauche de $A$ s'\'ecrivant $\Hom_D(E/V,E)$ pour un sous-$D$-espace vectoriel $V$ de $E$ alors
 $I^\circ/I^\circ I$ est canoniquement isomorphe \`a $\Hom_D(V,V)$. Si de plus, $J$ est un id\'eal \`a gauche de $A$ contenant $I$ s'\'ecrivant $\Hom_D(E/W,E)$, avec $W\subset V$ alors $I^\circ J/I^\circ I$ est canoniquement isomorphe \`a  $\Hom_D(V/W,V)$.
\end{lem}
\begin{proof}
Si $V\subset E$, on a la suite exacte :
$$
\xymatrix{0\ar[r]&V\ar[r]&E\ar[r]&E/V\ar[r]&0}
$$
et en appliquant le foncteur $\Hom_D(-,V)$ qui est exact, on obtient la suite exacte suivante :
$$
\xymatrix{0\ar[r]&\Hom_D(E/V,V)\ar[r]&\Hom_D(E,V)\ar[r]&\Hom_D(V,V)\ar[r]&0}.
$$
Puisque d'apr\`es \ref{dimension} $\Hom_D(E/V,V)=I^\circ I$, la suite exacte se r\'e\'ecrit :
$$
\xymatrix{0\ar[r]&I^\circ I\ar[r]&I^\circ \ar[r]&\Hom_D(V,V)\ar[r]&0}.
$$
Par suite : $I^\circ/I^\circ I$ est canoniquement isomorphe \`a $\Hom_D(V,V)$.\\
Supposons que $J=\Hom_D(E/W,E)$. Un raisonnement analogue appliqu\'e \`a la suite exacte courte suivante
$$
\xymatrix{0\ar[r]&V/W\ar[r]&E/W\ar[r]&E/V\ar[r]&0}
$$
 nous permet de voir que $I^\circ J/I^\circ I=\Hom_D(V/W,V)$.
\end{proof}

\begin{cor}\label{cor_quot_1}
Soit $A$ une alg\`ebre d'Azumaya sur $k$, de degr\'e $n$.
Si $I$ est un id\'eal \`a gauche de $A$, alors on a un isomorphisme canonique :
$$I^\circ/I^\circ I \to \End_A (I^\circ) $$
induite par l'application $x\mapsto (y\mapsto xy)$.
De plus : $\ind (\End_A I^\circ) = \ind A$ (\cite[proposition 1.10]{invol}). Le groupe $I^\circ /I^\circ I$ est muni d'une structure d'alg\`ebre d'Azumaya sur $k$ dont la classe est \'egale \`a celle de $A$ dans le groupe $\brauer (k)$ et est de degr\'e $n-i$ si $\dim_k I=ni$.
\end{cor}
\begin{proof}
On peut supposer que $A=\End_D(E)$ avec les notations pr\'ec\'edentes. $I$ s'\'ecrit alors $\Hom_D(E/V,E)$ pour un certain $V$.
D'apr\`es \cite[proposition 1.12]{invol}, $\Hom_D(V,V)$  s'identifie canoniquement \`a $\End_A (I^\circ)= \End_A (\Hom_D(E,V))$ gr\^ace \`a la multiplication \`a gauche. L'assertion sur le degr\'e est une cons\'equence du lemme \ref{dimension}.
\end{proof}

\begin{lem}\label{quot2}
Conservons les notations pr\'ec\'edentes.
Alors $J/J^\circ J$ est canoniquement isomorphe \`a $\Hom_D(E/W,E/W)$ et $I/J^\circ I$ est canoniquement isomorphe \`a $\Hom_D(E/V,E/W)$.
\end{lem}
\begin{proof}
La preuve est analogue \`a celle du lemme \ref{quot1}
\end{proof}
\begin{cor}\label{cor_quot_2}
Soit $A$ une alg\`ebre d'Azumaya sur $k$, de degr\'e $n$.
Si $I$ est un id\'eal \`a gauche de $A$, alors on a un isomorphisme canonique :
$$I/I^\circ I \to \End_A I$$
induite par l'application $x\mapsto (y\mapsto yx)$.
De plus : $\ind (\End_A I^\circ) = \ind A$. Le groupe $I/I^\circ I$ est muni d'une structure d'alg\`ebre d'Azumaya sur $k$ dont la classe est \'egale \`a celle de $A$ dans le groupe $\brauer (k)$ et est de degr\'e $i$ si $\dim_k I=ni$.
\end{cor}
\begin{proof}
La preuve est identique \`a celle du corollaire \ref{cor_quot_1}.
\end{proof}

\Subsection{Bijections}

La proposition suivante est la clef de la preuve du th\'eor\`eme \ref{thm_effectif} :

\begin{prop}\label{prop_effectif}
Soit $A$ une alg\`ebre d'Azumaya sur $k$. Notons $n$ son degr\'e. Soit $I\subset A$ un id\'eal \`a gauche de $k$-dimension $n$.\\
On a une bijection canonique~:
$$
\begin{array}{ccc}
\enscolquatre {J}
                               {I\subset J \subset A}
                               {J \mbox{ id\'eal \`a gauche}}
                               {\dim_k J=nj}
 &\rightleftarrows&
\enscolquatre{W}
                              {W \subset I^\circ I}
                             {W \mbox{ sous-espace vectoriel}}
                             {\dim_k W=(n-j)}
\\
J&\mapsto&J^\circ I
\end{array}
$$
\end{prop}
\begin{proof}
D'apr\`es le lemme \ref{dimension}, $\dim_k J^\circ I = (n-j)$, et donc on a bien une fl\`eche not\'ee $\phi$ :
$$
\xymatrix{\enscolquatre {J}
                               {I\subset J \subset A}
                               {J \mbox{ id\'eal \`a gauche}}
                               {\dim_k J=nj}
           \ar[r]^-\phi&
           \enscolquatre{W}
                              {W \subset I^\circ I}
                             {W \mbox{ sous-espace vectoriel}}
                             {\dim_k W=(n-j)} }
$$
Remarquons d\'ej\`a que l'existence de $I$ force $d=\ind A=1$. $A$ est donc d\'eploy\'ee. D'apr\`es le th\'eor\`eme de Wedderburn, il existe un $k$-espace vectoriel $E$ tel que l'on ait $A\simeq \End_k(E)$.
Par suite, gr\^ace au lemme \ref{id-vect}, on dispose de $V\subset E$ tel que $I$ s'identifie \`a $\Hom_k(E/V,E)$ via l'isomorphisme $A\simeq \End_k(E)$. En particulier $\dim_k V=n-1$. Puisque $\dim_k E/V=1$, on a une bijection canonique:
$$
\xymatrix{\enscol{U \subset \Hom(E/V,V)}{\mbox{sous-espace vectoriel de dimension $l$}} 
          \ar@<2pt>[r]^-f&
          \enscol{U \subset V}{\dim_k U=l}
          \ar@<2pt>[l]\\
          U\ar@{|->}[r]&\bigcup_{u\in U} \Im u}
$$
Par suite, le diagramme commutatif suivant permet de conclure :

$$
\xymatrix{
   \enscolquatre{J}
               {I\subset J \subset A}
               {J\mbox{ id\'eal \`a gauche}}
               {\dim_k J=nj}
   \ar[r]^-\phi
   \ar@<-2pt>[dd]_{\ref{id-vect}}
   &
   \enscolquatre{W}
                      {W \subset I^\circ I}
                      {W \mbox{ sous-espace vectoriel}}    
                      {\dim_k W =n-j}
   \ar@<-2pt>[d]_{\ref{id-vect}}
   \\
   &
   \enscolquatre{W}
                                 {W \subset \Hom(E/V,V)}
                                 {W\mbox{ sous-espace vectoriel}}
                                 {\dim_k W=n-j}
   \ar@<-2pt>[u]
   \ar@<-2pt>[d]_-{f}
   \\
   \enscolquatre{W}
                                 {W\subset V \subset E} 
                                 {W\mbox{ sous-espace vectoriel}}
                                 {\dim_k W=n-j}
   \ar@<-2pt>[uu]
   \ar@{=}[r]
   &   
   \enscolquatre{W}
                                 {W \subset V}
                                 {W\mbox{ sous-espace vectoriel}}
                                {\dim_k W=n-j}
  \ar@<-2pt>[u]}
$$
\end{proof}

\begin{rem}\label{rem_prop_effectif}
Notons $\phi$ la bijection de la proposition pr\'ec\'edente. On peut donner une expression explicte de $\phi^{-1}$. Avec les notations de la proposition, elle s'\'ecrit~:
$$
\begin{array}{ccc}
\enscolquatre{W}
                              {W \subset I^\circ I}
                             {W \mbox{ sous-espace vectoriel}}
                             {\dim_k W=(n-j)}
 &\rightleftarrows&
\enscolquatre {J}
                               {I\subset J \subset A}
                               {J \mbox{ id\'eal \`a gauche}}
                               {\dim_k J=nj}

\\
W&\mapsto&^\circ \! (WA).
\end{array}
$$
En effet, si on note par $\psi$ l'application ci-dessus, d'apr\`es le lemme \ref{rem_annulateurs}, le $k$-espace vectoriel $^\circ \! (WA)$ est un id\'eal \`a gauche de $A$.
Si $J\supset I$ est de $k$-dimension $nj$ alors $^\circ \! (J^\circ IA)= ^\circ \!(J^\circ)$, car par simplicit\'e de $A$, l'id\'eal $IA$ \'etant bilat\`ere et non vide, il est \'egal \`a $A$. Or $^\circ \!(J^\circ)=J$ gr\^ace au lemme \ref{lem_double_perp}, donc $\psi(\phi(J))=J$. Ainsi, puisque $\psi$ est un inverse \`a gauche de $\phi$ et puisque $\phi$ est surjective, alors pour tout $W \subset I^\circ I$, l'id\'eal \`a gauche $^\circ (WA)$ contient $I$ et est de $k$-dimension $n-j$, donc l'application $\psi$ est bien d\'efinie.

Soit $W$ un $k$-espace vectoriel de $I^\circ I$ de $k$-dimension $n-j$. On calcule  $\phi(\psi(W))=(^\circ \!(WA))^\circ I$.\
D'apr\`es le lemme \ref{lem_double_perp}, $\phi(\psi(W))=WAI=WI$. Or $WI \subset W$ et donc pour des raisons de dimension $\phi(\psi(W))=W$.

\end{rem}

\begin{lem}\label{bij_first}
Soit $A$ une alg\`ebre d'Azumaya sur $k$ de degr\'e $n$. Soient $r>1$ et $1 \leqslant i_1<\ldots<i_r \leqslant n$ des entiers. Soit $I_1\subset A$ un id\'eal de $k$-dimension $ni_1$ (i.e.: un point rationnel de $\SB{i_1}{A}$). On a alors une bijection canonique~:
$$
\xymatrix{
   \enscolquatre{(I_2,\ldots,I_r)}
                {I_1 \subset I_2 \subset \ldots \subset I_r \subset  A}      
                {I_j\mbox{ id\'eal \`a gauche}}
                {\dim_k I_j=ni_j}
   \ar@<2pt>[r]
   &
   \enscolquatre {(J_2,\ldots\,J_r)}
               {J_2\subset \ldots \subset J_r \subset I_1^\circ/I_1^\circ I_1}
               {J_j\mbox{ id\'eal \`a gauche}}
               {\dim_k J_j=(n-i_1)(i_j-i_1)}
   \ar@<2pt>[l]
   \\
   (I_2,\ldots,I_r)\ar@{|->}[r]
   &(I_1^\circ I_2 /I_1^\circ I_1,\ldots,I_1^\circ I_r /I_1^\circ I_1).}
$$
\end{lem}
\begin{proof}
Posons $B=\End_A (I_1^\circ)=I_1^\circ/I_1^\circ I_1$(corollaire \ref{cor_quot_1}). C'est une alg\`ebre d'Azumaya sur $k$ d'indice 
$d=\ind A$ et de degr\'e $n-i_1$.
Pour tout $j$, posons $J_j=I_1^\circ I_j/I_1^\circ I_1$. On a $J_j \subset B$, et les $J_j$ sont des id\'eaux \`a gauche de $B$ de $k$-dimension $(n-i_1)(i_j-i_1)$ d'apr\`es \ref{dimension}. 
On a alors une fl\`eche $\phi$~: 
$$
\xymatrix{
  \enscolquatre{(I_2,\ldots , I_r)}
                    {I_1 \subset I_2 \subset \ldots \subset I_r \subset A}
                  {I_j\mbox{ id\'eal \`a gauche}}
                 {\dim_k I_j=ni_j}
  \ar[r]^-\phi
  &
  \enscolquatre {(J_2,\ldots,J_r)}
                     {J_2\subset \ldots \subset J_r \subset I_1^\circ/I_1^\circ I_1}
                  {J_j\mbox{ id\'eal \`a gauche}}
                  {\dim_k J_j=(n-i_1)(i_j-i_1)} }
$$
On peut supposer que $A=\End_D(E)$ pour un corps gauche $D$ de dimension $d^2$ sur $k$ et un $D$-espace vectoriel \`a droite $E$ de dimension $r$. Alors $n=rd$. Soit $V_1\subset E$ tel que $I_1=\Hom_D(E/V_1,E)$. Ainsi $B=\End_A (I_1^\circ)=I_1^\circ/I_1^\circ I_1=\Hom(V_1,V_1)$ (lemme \ref{quot1}). 

Le diagramme commutatif suivant dont les fl\`eches verticales sont donn\'ees par la proposition \ref{id-vect} permet de conclure :
$$
\xymatrix{
   \enscolquatre{(I_2,\ldots,I_r)}
                                 {I_1 \subset I_2 \subset \ldots \subset I_r \subset A}
                                 {I_j\mbox{ id\'eal \`a gauche}}
                                 {\dim_k I_j=ni_j}
   \ar[r]^-\phi
   \ar@<-2pt>[d]
   &
   \enscolquatre{(J_2,\ldots,J_r)}
                                 {J_2\subset \ldots \subset J_r \subset I_1^\circ/I_1^\circ I_1}
	                         {J_j\mbox{ id\'eal \`a gauche}}
                                 {\dim_k J_j=(n-i_1)(i_j-i_1)}                    
   \ar@<-2pt>[d]
   \\
   \enscolquatre {(V_2,\ldots, V_r)}
	                          {E\supset V_1\supset V_2\supset\ldots \supset V_r}
                                  {V_j \mbox{ sous-espace vectoriel}}
                                  {\dim_D V_j=\frac{n-i_j}{d}}
   \ar@{=}[r]
   \ar@<-2pt>[u]
   &
   \enscolquatre {(W_2, \ldots , W_r)}
	                          {V_1 \supset W_2 \supset \ldots \supset W_r}
	                          {W_j \mbox{ sous-espace vectoriel}}
                                  {\dim_D W_j=\frac{(n-i_1)-(i_j-i_1)}{d}}
   \ar@<-2pt>[u] }
$$
\end{proof}
\begin{rem}\label{rem_bij_first}
Notons $\phi$ la bijection du lemme pr\'ec\'edent. On peut d\'ecrire de mani\`ere explicite sa bijection r\'ecriproque. Avec les notations du lemme, c'est l'application $\psi$ suivante~:
$$
\xymatrix{
   \enscolquatre{(J_2,\ldots,J_r)}
                {J_2 \subset \ldots \subset J_r \subset  I_1^\circ/I_1^\circ I_1}      
                {J_j\mbox{ id\'eal \`a gauche}}
                {\dim_k J_j=(n-i_1)(i_j-i_1)}
   \ar@<2pt>[r]
   &
   \enscolquatre {(I_2,\ldots\,I_r)}
               {I_1 \subset I_2\subset \ldots \subset I_r \subset A}
               {I_j\mbox{ id\'eal \`a gauche}}
               {\dim_k I_j=ni_j}
   \ar@<2pt>[l]
   \\
   (J_2,\ldots,J_r)\ar@{|->}[r]
   &(A\pi_1^{-1}(J_2),\ldots,A\pi_1^{-1}(J_r))}
$$
o\`u $\pi_1$ est la projection $A\to A/I_1^\circ I_1$.
Tout d'abord pour tout $j$, $A\pi_1^{-1}(J_j)$ est bien un id\'eal \`a gauche de $A$, et l'on a bien $A\pi_1^{-1}(J_2)\subset \ldots \subset A\pi_1^{-1}(J_r)$. Fixons un $j\in \{2,\ldots, r\}$, si $I_j$ est un id\'eal \`a gauche de $A$ de $k$-dimension $ni_j$, alors $A\pi_1^{-1}(I_1^\circ I_j/I_1^\circ I_1) = AI_1^\circ I_j=I_j$, car, par simplicit\'e de $A$, $AI_1^\circ=A$. Ainsi d'une part $\psi$ est un inverse \`a gauche de $\phi$ et d'autre part puisque $\phi$ est surjective, l'id\'eal $A\pi_1^{-1}(J_j)$ est de $k$-dimension $ni_j$ d\`es que $J_j$ est de $k$-dimension $(n-i_1)(i_j-i_1)$ et contient $I_1$.

Soient  $J_2,\ldots,J_r$ des id\'eaux \`a gauche de $I_1^\circ/I_1^\circ I_1$ de $k$-dimension respective $(n-i_1)(i_j-i_1)$ tels que $J_2 \subset \ldots \subset J_r$. Pour tout $j\in \{2\ldots r\}$, on a~: $I_1^\circ (A\pi_1^{-1}(J_j)) /I_1^\circ I_1=\pi^{-1}_1(J_j)/I_1^\circ I_1$ puisque $I_1^\circ A=A$; et $\pi^{-1}_1(J_j)/I_1^\circ I_1 = J_j$. Ainsi $\phi(\psi(J_2,\ldots,J_r))=(J_2,\ldots,J_r)$.
\end{rem}

De m\^eme, on a :
\begin{lem}\label{bij_last}
Soit $A$ une alg\`ebre d'Azumaya sur $k$ de degr\'e $n$. Soient $r>1$ et $1 \leqslant i_1<\ldots<i_r\leqslant n$ des entiers. Soit $I_r\subset A$ un id\'eal de $k$-dimension $ni_r$ (i.e.: un point rationnel de $\SB{i_r}{A}$). On a alors une bijection canonique~:
$$
\xymatrix{
   \enscolquatre{(I_1,\ldots , I_{r-1})}
                      {I_1 \subset \ldots \subset I_{r-1} \subset I_r \subset  A}      
                      {I_j\mbox{ id\'eal \`a gauche}}
                      {\dim_k I_j=ni_j}
   \ar@<2pt>[r]
   &
   \enscolquatre {(J_1,\ldots,J_{r-1})}
                    {J_1\subset \ldots \subset J_{r-1} \subset I_r/I_r^\circ I_r}
                    {J_j\mbox{ id\'eal \`a gauche}}
                  {\dim_k J_j=i_ji_r}
   \ar@<2pt>[l]
   \\
   (I_j)_{j=1\ldots (r-1)}
    \ar@{|->}[r]
   & 
   (I_j/I_r^\circ I_j)_{j=1\ldots (r-1)} }
$$
\end{lem}

\begin{proof}
La preuve est analogue \`a celle qui pr\'ec\`ede.\\
Posons $B=\End_A (I_r)=I_r/I_r^\circ I_r$ (corollaire \ref{cor_quot_2}). C'est donc une alg\`ebre d'Azumaya sur $k$ d'indice $d=\ind A$ et de degr\'e $i_r$. Pour $j\leqslant (r-1)$, posons $J_j=I_j/I_r^\circ I_j$. La multiplication nous donne une action de $I_r$ (\`a gauche) sur $I_j$. Cette action donne naissance \`a une action de $I_r$ sur $J_j$ puisque $I_r I_r^\circ=0$. Finalement, elle permet d'obtenir une action (multiplication \`a gauche) de $B$ sur $J_j$. Ainsi $J_j$ est muni d'une structure de $B$-module. De plus, l'injection canonique $I_j \to I_r$ nous donne une application $I_j \to I_r/I_r^\circ I_r$ dont le noyau est $I_r^\circ I_j$, ainsi on a une application injective  $J_j \to B$. Cette injection est compatible avec la structure de $B$-module et nous appelerons encore $J_j$ son image dans $B$. Ainsi, les $J_j$ sont des id\'eaux \`a gauche de $B$ de $k$-dimension $i_ji_r$ d'apr\`es \ref{dimension}. On a donc une fl\`eche $\phi$~:

$$
\xymatrix{
   \enscolquatre{(I_1,\ldots,I_{r-1})}
                                {I_1 \subset \ldots \subset I_{r-1} \subset I_r \subset  A}      
                                {I_j\mbox{ id\'eal \`a gauche}}
                                {\dim_k I_j=ni_j}
   \ar[r]^-\phi
   &
   \enscolquatre{(J_1,\ldots,J_{r-1})}
                                 {J_1\subset \ldots \subset J_{r-1} \subset I_r/I_r^\circ I_r}
                                 {J_j\mbox{ id\'eal \`a gauche}}
                                 {\dim_k J_j=i_ji_r}
   }
$$

On peut supposer que $A=\End_D(E)$ pour un corps gauche $D$ de dimension $d^2$ sur $k$ et un $D$-vectoriel \`a droite $E$ de dimension $r$. Alors $n=rd$. Soit $V_r\subset E$ tel que $I_r=\Hom_D(E/V_r,E)$. On a alors $B=\End_A (I_r)=I_r/I_r^\circ I_r=\Hom(E/V_r,E/V_r)$ (corollaire \ref{cor_quot_2}).

Le diagramme commutatif suivant dont les fl\`eches verticales sont donn\'ees par la proposition \ref{id-vect} permet de conclure :
$$
\xymatrix{
   \enscolquatre{(I_2,\ldots, I_{r-1})}
                      {I_1 \subset \ldots \subset I_{r-1}\subset I_r \subset A}
                      {I_j\mbox{ id\'eal \`a gauche}}
                      {\dim_k I_j=ni_j}
   \ar[r]^-\phi
   \ar@<-2pt>[d]
   &
  \enscolquatre{(J_1,\ldots,J_{r-1})}
                                {J_1\subset \ldots \subset J_{r-1} \subset I_r/I_r^\circ I_r}
                                {J_j\mbox{ id\'eal \`a gauche}}
                                {\dim_k J_j=i_ji_r}
   \ar@<-2pt>[d]
   \\
   \enscolquatre{(V_1, \ldots , V_{r-1})}
                                 {E\supset V_1\supset V_2\supset\ldots \supset V_r}
                                 {V_j \mbox{ sous-espace vectoriel}}
                                 {\dim_D V_j=\frac{n-i_j}{d}}
   \ar@<2pt>[r]^-{\pi}
   \ar@<-2pt>[u]
   &
   \enscolquatre{(W_1 , \ldots , W_{r-1})}
                                 {E/V_r \supset W_1 \supset \ldots \supset W_{r-1}}
                                 {W_j \mbox{ sous-espace vectoriel}}
                                 {\dim_D W_j=\frac{i_r-i_j}{d}}
   \ar@<-2pt>[u]
   \ar@<2pt>[l] 
   }
$$
o\`u $\pi$ est la projection naturelle : $E\to E/V_r$.
\end{proof}

\begin{rem}\label{rem_bij_last}
Ici aussi il est possible de donner une expression explicite de la bijection r\'eciproque. Conservons les notations du lemme, La bijection r\'eciproque s'\'ecrit~:
$$
\xymatrix{
   \enscolquatre{(J_1,\ldots , J_{r-1})}
                      {J_1 \subset \ldots \subset J_{r-1} \subset  I_r/I_r^\circ I_r }      
                      {J_j\mbox{ id\'eal \`a gauche}}
                      {\dim_k J_j=i_ji_r}
   \ar@<2pt>[r]
   &
   \enscolquatre {(I_1,\ldots,I_{r-1})}
                    {I_1\subset \ldots \subset I_{r-1} \subset I_r \subset A}
                    {I_j\mbox{ id\'eal \`a gauche}}
                  {\dim_k I_j=ni_j}
   \ar@<2pt>[l]
   \\
   (J_j)_{j=1\ldots (r-1)}
    \ar@{|->}[r]
   & 
   (A\pi_r^{-1}(J_j))_{j=1\ldots (r-1)} }
$$
o\`u $\pi_r$ est la projection canonique $A\to A/(I_r^\circ I_r)$.
La preuve est analogue \`a celle donn\'ee en \ref{rem_bij_first}.
\end{rem}

La proposition suivante est la clef de la preuve du th\'eor\`eme \ref{thm_general}.
\begin{prop}\label{bij_general}
Soit $A$ une alg\`ebre d'Azumaya sur $k$ de degr\'e $n$. Soient $r>1$ et $i_1<\ldots<i_r\leqslant n$ des entiers. Soit $s\in\lbrace 1\ldots r \rbrace$. Soit $I_s\subset A$ un id\'eal de $k$-dimension $ni_s$ (i.e.: un point rationel de $\SB{i_s}{A}$). L'application~:
$$
\xymatrix@-1pc{
    (I_1,\ldots,\hat{I_s},\ldots,I_r)
    \ar@{|->}[d]
    \\
   ((I_j/I_r^\circ I_j)_{j=1..(s-1)},
    (I_s^\circ I_j /I_s^\circ I_s)_{j=(s+1)..r} )}
$$
d\'efinit une bijection~:
$$
 \enscolquatre{(I_1,\ldots, \widehat{I_s},\ldots, I_r)}
                                 {I_1 \subset \ldots \subset I_r \subset A}
                                 {I_j \mbox{ id\'eal \`a gauche}}
                                 {\dim_k I_j=ni_j}
\leftrightarrows
\begin{array}{c}
   \enscolquatre{(J_1,\ldots,J_{s-1})}
                                 {J_1\subset \ldots \subset J_{s-1} \subset I_s/I_s^\circ I_s}
                                 {J_j \mbox{ id\'eal \`a gauche}}
                                {\dim_k J_j=i_ji_s} \\
\times \\
 \enscolquatre {(J_{s+1},\ldots,J_r)}
                                  {J_{s+1}\subset \ldots\subset J_r \subset I_s^\circ/I_s^\circ I_s}
                                  {J_j \mbox{ id\'eal \`a gauche}}
                                 {\dim_k J_j=(n-i_s)(i_j-i_s)} 
\end{array}
$$
%
%$$
%\xymatrix{
%   \enscolquatre{(I_1,\ldots, \widehat{I_s},\ldots, I_r)}
%                                 {I_1 \subset \ldots \subset I_r \subset A}
%                                 {I_j \mbox{ id\'eal \`a gauche}}
%                                 {\dim_k I_j=ni_j}
%\ar@<2pt>[d]
%   \\
%   \enscolquatre{(J_1,\ldots,J_{s-1})}
%                                 {J_1\subset \ldots \subset J_{s-1} \subset I_s/I_s^\circ I_s}
%                                 {J_j \mbox{ id\'eal \`a gauche}}
%                                {\dim_k J_j=i_ji_s}
%   \times
%   \enscolquatre {(J_{s+1},\ldots,J_r)}
%                                  {J_{s+1}\subset \ldots\subset J_r \subset I_s^\circ/I_s^\circ I_s}
%                                  {J_j \mbox{ id\'eal \`a gauche}}
%                                 {\dim_k J_j=(n-i_s)(i_j-i_s)}
%\ar@<2pt>[u]}
%$$
\end{prop}
\begin{proof}
Ce sont les deux lemmes pr\'ec\'edents.
\end{proof}

\begin{rem}\label{rem_bij}
Les propositions et lemmes  pr\'ec\'edents s'\'etendent au cas o\`u $A$ est une alg\`ebre d'Azumaya sur une $k$-alg\`ebre $R$. En effet, les fl\`eches $\phi$ introduites \`a la proposition \ref{prop_effectif} puis aux lemmes \ref{bij_first} et \ref{bij_last} et enfin \`a la proposition \ref{bij_general} ont toujours un sens si $A$ est une alg\`ebre d'Azumaya sur une $k$-alg\`ebre $R$. 
Il en est de m\^eme de leur r\'eciproque, donn\'ees dans les remarques \ref{rem_prop_effectif}, \ref{rem_bij_first} et \ref{rem_bij_last}.
%Puis, pour montrer que ce sont encore des bijections, il suffit de se ramener par %localisation au cas o\`u $R$ est un corps. Dans ce cas, $A$ est une alg\`ebre %d'Azumaya sur un corps et les r\'esultats de la section montrent que les fl\`eches %$\phi$ sont encore des bijections.
\end{rem}
\Subsection{Preuve du th\'eor\`eme \ref{thm_effectif}}

Gr\^ace aux r\'esultats du paragraphe pr\'ec\'edent, on est en mesure de prouver le  th\'eor\`eme \ref{thm_effectif}:

\begin{proof}
Fixons $1=i_1< i_2<\ldots<i_r \leqslant n$ une suite d'entiers. Notons $S=\SB{i_1}{A}=\SB{}{A}$ et $X=\Flag{1,i_2,\ldots,i_r}{A}$. Et posons, pour $j>1$, $a_j=(n-i_j)$.
\noindent
Consid\'erons l'application identit\'e : $S \to S$. Elle nous fournit de facto un $S$-point de $S$ i.e. un \'el\'ement de $S(S)$. Il correspond \`a un faisceau d'id\'eaux $\mathcal{I} \subset A\otimes_k \mathcal{O}_S$. Soit $s$ un $L$-point de $S$, le diagramme suivant est commutatif 
$$
\xymatrix{
   \Spec{L}
   \ar[r]^-{s}
   \ar@/_1pc/[rr]_-{s}
   &
   S
   \ar[r]^-{\mathcal{I}}
   &
   S}.
$$
Or $s$ correspond \`a un id\'eal \`a gauche de $k$-dimension $n$, $I_s$ de $A\otimes_k L$, on a : $\mathcal{I}_s=I_s$. Ou de mani\`ere abusive : $\mathcal{I}_I=I$ en confondant $I=I_s$ et $s$, abus que l'on fera dans la suite si aucune ambigu\"it\'e n'est \`a craindre.
Comme nous avons d\'efini l'annulateur de $I$, $I^\circ$, nous pouvons d\'efinir l'annulateur de $\mathcal{I}$, $\mathcal{I}^\circ$. On consid\`ere alors le fibr\'e sur $S$ suivant : $\mathcal{V}=\mathcal{I}^\circ \mathcal{I}$. Pour tout point $I$ de $S$, on a :
$\mathcal{V}_I=I^\circ I$.
C'est un fibr\'e vectoriel de rang $n-1$ d'apr\`es le lemme 
\ref{dimension}.\\

Consid\'erons la fl\`eche suivante :
$$
(I_1,\ldots\,I_r)\to(I_r^\circ I_1,I_{(r-1)}^\circ I_1,\ldots,I_2^\circ I_1\subset I_1^\circ I_1)
$$
Elle permet de d\'efinir pour toute $k$-alg\`ebre $R$ une application bijective (proposition \ref{prop_effectif} et remarque \ref{rem_bij}) $X(R) \to \Drapeau{n-i_r,\ldots,n-i_2}{\mathcal{V}}(R)$. Et par suite, elle correspond \`a un isomorphisme~: $$X\to \Drapeau{n-i_r,\ldots,n-i_2}{\mathcal{V}}.\qed$$\noqed
\end{proof}

\Subsection{Preuve du th\'eor\`eme \ref{thm_general}}

\begin{proof}
Fixons $1 \leqslant i_1< i_2<\ldots<i_r \leqslant n$ une suite d'entiers. Notons $S=\SB{i_s}{A}$ et $X=\Flag{i_1,i_2,\ldots,i_r}{A}$.

\`A partir du fibr\'e $\mathcal{I}$ sur $S$ d\'efini dans la preuve du th\'eor\`eme \ref{thm_general}, d\'efinissons $\mathcal{B}_{+}=\mathcal{I}^\circ/\mathcal{I}^\circ\mathcal{I}$ 
et $\mathcal{B}_{-}=\mathcal{I}/\mathcal{I}^\circ\mathcal{I}$. 
D'apr\`es le lemme \ref{quot1} et le corollaire \ref{cor_quot_1}, $\mathcal{B}_{+}$ est une alg\`ebre d'Azumaya sur $S$ de degr\'e $n-i_s$  canoniquement  isomorphe \`a $\End_S \mathcal{I}^\circ$. De m\^eme (\ref{quot2} et \ref{cor_quot_2}), $\mathcal{B}_{-}$ est une alg\`ebre d'Azumaya sur $S$ de degr\'e $i_s$  canoniquement  isomorphe \`a $\End_S \mathcal{I}$.

Consid\'erons la fl\`eche suivante :
$$
\xymatrix{
    (I_1,\ldots,\hat{I_s},\ldots,I_r)
    \ar@{|->}[d]
    \\
   ((I_j/I_r^\circ I_j)_{j=1..(s-1)},
    (I_s^\circ I_j /I_s^\circ I_s)_{j=(s+1)..r} )}
$$
Elle permet de d\'efinir pour toute $k$-alg\`ebre $R$ une application bijective (proposition \ref{bij_general} et remarque \ref{rem_bij}) :

$$X(R) \to
\Flag{i_1,\ldots,i_{s-1}}{\mathcal{B_{-}}}
  \times_{\SB{i_s}{A}}
  \Flag{i_{s+1}-i_s,\ldots,i_r-i_s}{\mathcal{B_{+}}}(R)
.$$
Et par suite, elle correspond \`a un isomorphisme

$$X \to
\Flag{i_1,\ldots,i_{s-1}}{\mathcal{B_{-}}}
  \times_{\SB{i_s}{A}}
  \Flag{i_{s+1}-i_s,\ldots,i_r-i_s}{\mathcal{B_{+}}}
.\qed$$\noqed
\end{proof}

\Subsection{Preuve du corollaire \ref{cor_thm_general}}

%\begin{lem}
%Soient $S$ un $k$-sch\'ema, $X\to S$ et $Y \to S$ deux fibrations en vari\'et\'e de drapeaux, alors $X\times_S Y \to S$ est une fibration en produits de vari\'et\'es de drapeaux.
%\end{lem}
%\begin{proof}
%Par d\'efinition-m\^eme du produit, on a le carr\'e cart\'esien suivant : 
%$$
%\xymatrix{
%  X\otimes_S Y
%  \ar[r]^{\pi_X}
%  \ar[d]^{\pi_Y}
%  \cartesien
%  &
%  X
 % \ar[d]^-{f_X}
 % \\
 % Y
 % \ar[r]^-{f_Y}
  %&
 % S}
%$$
%Puisque $f_Y\,:\,Y \to S$ est une fibration en vari\'et\'e de drapeaux, alors il en est de m\^eme pour $\pi_X\,:%\,X\otimes_S Y \to X$. Ainsi est-on ramen\'e \`a montrer : 
%Si $Y\to X$ et $X\to S$ sont des fibrations en vari\'et\'es de drapeaux, alors la compos\'ee $Y\to S$ est une %fibration en produit de vari\'et\'es de drapeaux, ce qui est clair.
%\end{proof}

On est en mesure de prouver le corollaire \ref{cor_thm_general} : 
\begin{proof}
Ici \brauer(S) d\'esigne le groupe de Brauer de $S$, puisque $S$ est lisse sur $k$, $S$ est un sch\'ema r\'egulier. Alors, d'apr\`es \cite[corollaire 1.10]{groth_dix_2}la fl\`eche naturelle d'\'evaluation
$$
\brauer (S) \to \brauer k(S)
$$ 
est injective, o\`u $k(S)$ d\'esigne le corps des fonctions de $S$. Or pour tout $I$, $\mathcal{B}_{+\,I}=\End_A I^\circ$ et 
$\mathcal{B}_{-\,I}=\End_A I$. Ainsi, puisque $\ind  A_{\kappa(I)}$=1 et que d'apr\`es le corollaire \ref{cor_quot_1}, les classes de $A_{\kappa(I)}$ et de $\End_A I^\circ$ sont \'egales dans 
$\brauer \kappa(I)$,  $\mathcal{B}_{+}$ est triviale dans $\brauer S$. De m\^eme, $\mathcal{B}_{-}$ est triviale dans $\brauer S$.
Ainsi les deux fl\`eches suivantes 

$$
\Flag{i_1,\ldots,i_{s-1}}{\mathcal{B^{-}}}
\to
\SB{i_s}{A}
$$
et
$$\Flag{i_{s+1}-i_s,\ldots,i_r-i_s}{\mathcal{B^{+}}}
\to
\SB{i_s}{A}
$$
sont des fibrations en vari\'et\'es de drapeaux non tordues, et par suite
$$
\xymatrix{
\Flag{i_1,\ldots,i_{s-1}}{\mathcal{B_{-}}}
  \times_{\SB{i_s}{A}}
  \Flag{i_{s+1}-i_s,\ldots,i_r-i_s}{\mathcal{B_{+}}}
\ar[d]
\\
\SB{i_s}{A} }
$$
est une fibration en produit de vari\'et\'es de drapeaux non tordues.
\end{proof}

\section{Calcul des groupes de Chow de fibr\'es en drapeaux}

\Subsection{Partitions d'entiers.}
\begin{defi}\label{defi_partition}
Si $n$ est un entier $\geqslant 1$, alors une partition de $n$ (cf. \cite[Paragraphe 2.1, d\'efinition A]{comtet_comb}) est la repr\'esentation de $n$ comme somme d'entiers $\geqslant 1$, sans consid\'eration d'ordre dans la somme. Nous noterons par $p(n)$ le nombre de telles repr\'esentations de $n$. L'entier $P(n,m)$ sera le nombre de telles repr\'esentations avec exactement $m$ termes dans la somme (on parlera alors de $m$-partition stricte de $n$), et $p(n,m)$ le nombre de telles repr\'esentations avec au plus $m$ termes dans la somme, on parlera alors de $m$-partition de $n$ :
$$p(n,m)=\sum_{i=1}^{m}P(n,i).$$
Si de plus $A$ est un entier, une $(m,A)$-partition (resp. stricte) de $n$ est une $m$-partition (resp. stricte) de $n$ telle que chaque terme intervenant dans la somme soit inf\'erieur \`a $A$. Nous noterons $P(n,m,A)$ le nombre de telles partitions strictes , et $p(n,m,A)$ le nombre de $(m,A)$-partition de $n$ : 
$$p(n,m,A)=\sum_{i=1}^{m}P(n,i,A).$$
\end{defi}
\begin{rem}
si $A\geqslant n$ alors $P(n,m,A)=P(n,m).$

Pour des entiers $n,m,A$, on a : $$p(n,m,A)=\sum_{i=1}^{m} P(n,i,A).$$ Ainsi : 
$$P(n,m,A)=p(n,m,A)-p(n,m-1,A).$$
\end{rem}

On dispose de plusieurs caract\'erisations pour les $m$-partitions  (\cite[d\'efinition B et th\'eor\`eme A paragraphe 2.1]{comtet_comb}) :
\begin{rem}
Soit $n$ un entier $\geqslant 1$. Se donner une partition de $n$ en exactement $m$ termes (au sens de \ref{defi_partition}) est \'equivalent \`a se donner une solution de
$$
\left\lbrace
\begin{array}{rcl}
y_1\geqslant\ldots\geqslant y_m &\geqslant& 1,\\
y_1+\ldots+y_m&=&n.
\end{array}
\right.
$$
Et cela est encore \'equivalent \`a se donner une solution de :
$$
\left\lbrace
\begin{array}{rcl}
x_1,\ldots,x_n &\geqslant& 0,\\
x_1+2x_2+\ldots+nx_n&=&n,\\
x_1+x_2+\ldots+x_n&=&m.
\end{array}
\right.
$$
\end{rem}
On a des caract\'erisations semblables pour les $(m,A)$-partitions :

\begin{rem}\label{rem_partitionA_x}
Soient $n,m,A\geqslant 1$ des entiers. Se donner une $(m,A)$-partition stricte de $n$ est \'equivalent \`a se donner une solution du probl\`eme suivant : 
$$
\left\lbrace
\begin{array}{rcl}
A\geqslant y_1\geqslant\ldots\geqslant y_m &\geqslant& 1,\\
y_1+\ldots+y_m&=&n.
\end{array}
\right.
$$
Ceci est encore \'equivalent \`a se donner une solution du probl\`eme suivant : 
$$
\left\lbrace
\begin{array}{rcl}
x_1,\ldots,x_n &\geqslant &0,\\
x_1+2x_2+\ldots+Ax_A&=&n,\\
x_1+x_2+\ldots+x_A&=&m.
\end{array}
\right.
$$
\end{rem}

\begin{lem}
Soient $n,m,A \geqslant 1$ des entiers.
\begin{enumerate}
\item[(i)]
$P(n,m,1)=1$ si $m=n$ et $0$ sinon,\\
$p(n,m,1)=1$ si $n\geqslant m$ et $0$ sinon.

\item[(ii)]
$P(n,1,A)=1$ si $A\geqslant n$,\\
$p(n,1,A)=1$ si $A\geqslant n$.

\item[(iii)]
$$
P(n,m,A) \geqslant 1 \iff  m \leqslant n\leqslant mA
$$
et
$$
p(n,m,A) \geqslant 1 \iff  n\leqslant mA
.$$
\end{enumerate}
\end{lem}
\begin{proof}
Supposons $n,m,A \geqslant 1$ donn\'es.
Les assertions (i) et (ii) sont \'evidentes.
Pour (iii), $P(n,m,A)\geqslant 1$ signifie que l'on dispose de $m$ entiers compris entre $1$ et $A$, disons $x_1,\ldots, x_m$
tels que $n=x_1+\ldots+x_m$. N\'ecessairement, on doit avoir $m \leqslant n\leqslant mA$. Inversement, si l'on suppose que $m \leqslant n\leqslant mA$, il est imm\'ediat de voir que $P(n,m,A) \geqslant 1$.
Pour $p(n,m,A)$, il suffit de remarquer que $p(n,m,A)=\sum_{i=1}^m P(n,i,A)$.
\end{proof}

Le th\'eor\`eme suivant est une g\'en\'eralisation du th\'eor\`eme B du paragraphe 2.1 de \cite{comtet_comb}.

\begin{thm}
Soient $n,m,A \geqslant 1$ des entiers tels que $n>m$, $m>1$ et $A>1$.
Alors on a : 
$$
p(n,m,A)=p(n,m-1,A)+p(n-m,m,A-1)
$$
\end{thm}
\begin{proof}
L'entier $p(n,m,A)$ est le nombre de solutions du probl\`eme suivant : 
$$
\left\lbrace
\begin{array}{rcl}
x_1,\ldots,x_n& \geqslant &0,\\
x_1+2x_2+\ldots+Ax_A&=&n,\\
x_1+x_2+\ldots+x_A&\leqslant &m.
\end{array}
\right.
$$
\noindent
Or ces solutions se d\'ecomposent en 2 groupes, celui consistant en les solutions telles que  $x_1+\ldots+x_A =m$ et celui consistant en les solutions telles que $x_1+\ldots+x_A < m$. Dans le second groupe, il y a exactement $p(n,m-1,A)$ \'el\'ements. Il reste donc \`a conna\^itre le cardinal du premier groupe. Or se donner une solution de : 
$$
\left\lbrace
\begin{array}{rcl}
x_1,\ldots,x_A &\geqslant &0,\\
x_1+2x_2+\ldots+Ax_A&=&n,\\
x_1+x_2+\ldots+x_A&=&m
\end{array}
\right.
$$
est \'equivalent \`a se donner une solution de : 
$$
\left\lbrace
\begin{array}{rcl}
A\geqslant y_1 \geqslant \ldots \geqslant y_m &\geqslant& 1,\\
y_1+\ldots+y_m&=&n,
\end{array}
\right.
$$
qui est \'equivalent \`a se donner une solution pour $s\leqslant m$ de :
$$
\left\lbrace
\begin{array}{rcl}
A-1\geqslant y_1 \geqslant \ldots \geqslant y_s &\geqslant& 1,\\
y_1+\ldots+y_s&=&n-m.
\end{array}
\right.
$$
Donc le cardinal du premier groupe est $p(n-m,m,A-1)$.
\end{proof}

\begin{cor}{\cite[Paragraphe 2.1, Th\'eor\`eme B]{comtet_comb}}
Pour des entiers $n,m\geqslant 1$ tels que $n>m$ et $m>1$, on a : 
$$p(n,m)=p(n,m-1)+p(n-m,m).$$
\end{cor}
\begin{proof}
Il suffit de prendre $A$ assez grand dans le th\'eor\`eme pr\'ec\'edent.
\end{proof}

Le th\'eor\`eme suivant est une variante du th\'eor\`eme d'Euler calculant la fonction g\'en\'eratrice de $P(n)$. Il est \`a rapprocher du th\'eor\`eme A du paragraphe 2.2 de \cite{comtet_comb}.
\begin{thm}
Fixons un entier $A\in \mathbf{N}^*$. La fonction g\'en\'eratrice de $P(n,m,A)$ (appartenant \`a $\mathbf{Z}\lbrack\lbrack x,y \rbrack
\rbrack $) est :
$$
\Xi_A(x,y)=\prod_{i=1}^{A} \frac{1}{1-xy^i} =\sum_{n,m}P(n,m,A)x^my^n.
$$
\end{thm}
\begin{proof}
D\'eveloppons le produit : 
\begin{eqnarray*}
\prod_{i=1}^{A}\frac{1}{1-xy^i}&=&\prod_{i=1}^{A} (\sum_{j_i\geqslant 0} x^{j_i}y^{ij_i} )\\
&=&\sum_{x_1,\ldots,x_A} x^{x_1+\ldots+x_A}y^{x_1+2x_2+\ldots+Ax_A}\\
\end{eqnarray*}
Il suffit alors de consid\'erer les diff\'erentes caract\'erisations des partitions, voir la remarque \ref{rem_partitionA_x}.
\end{proof}

\begin{defi}
Soit $n\geqslant 1$ un entier. Soient $r\geqslant 1$, $m_1\ldots, m_r$ et $A_1,\ldots,A_r$ $2r$ entiers, une $((m_1,A_1),\ldots,(m_r,A_r))$-partition stricte de $n$ est une solution du syst\`eme en les $y_{i,j}$ suivant: 
$$
\left\lbrace
\begin{array}{l}
   \forall i =1\ldots r
   \left\lbrace
   \begin{array}{rcl}
   A_i\geqslant y_{i,1} \geqslant \ldots \geqslant y_{i,m_i} &\geqslant& 1\\
   y_{i,1}+\ldots+y_{i,m_i}&=&n_i
   \end{array}
   \right.
  \\
  n_1,\ldots,n_r\geqslant 0\\
  n_1+\ldots+n_r=n
\end{array}
\right.
$$
Nous noterons $Q(n,(m_1,A_1),\ldots,(m_r,A_r))$ le nombre de telles partitions.
L'entier $q(n,(m_1,A_1),\ldots,(m_r,A_r))$ sera le cardinal des $(s_1,A_1),\ldots,(s_r,A_r)$-partitions strictes de $n$ pour $s_1\leqslant m_1,\ldots,s_r\leqslant m_r$, i.e. le nombre de solutions du probl\`eme suivant : 
$$
\left\lbrace
\begin{array}{l}
   \forall i =1\ldots r
   \left\lbrace
   \begin{array}{rcl}
   A_i\geqslant y_{i,1} \geqslant \ldots \geqslant y_{i,m_i} &\geqslant& 0\\
   y_{i,1}+\ldots+y_{i,m_i}&=&n_i
   \end{array}
   \right.
  \\
  n_1,\ldots,n_r\geqslant 0\\
  n_1+\ldots+n_r=n
\end{array}
\right.
$$
On parlera alors de $((m_1,A_1),\ldots,(m_r,A_r))$-partition de $n$.
\end{defi}

\begin{nota}
Pour des entiers $m,A$, on notera $\pi_{(m,A)}$ la distribution associ\'ee \`a la fonction $p$ : 
$$
\begin{array}{rcl}
\mathbf{N}&\to&\mathbf{N}\\
i&\mapsto&\pi_{(m,A)}(i)=p(i,m,A).
\end{array}
$$
\end{nota}

Rappelons la d\'efinition du produit de convolution : 
\begin{defi}
Soient $f,g$ deux fonctions sur $\mathbf{N}$. Le produit de convolution est la fonction sur $\mathbf{N}$ d\'efinie de la fa\c{c}on suivante :
$$
\begin{array}{rcl}
\mathbf{N}&\to&\mathbf{N}\\
n&\mapsto&(f\star g) (n)=\sum_{l=0}^{n} f(n-l)g(l).
\end{array}
$$
\end{defi}

\begin{rem}
Reprenons les notations de la d\'efinition pr\'ec\'edente. Si l'on d\'efinit les deux s\'eries formelles (appartenant \`a l'anneau $\mathbf{Z}[[t]]$) :
$$
\begin{array}{ccc}
F(t)&=&\sum_n f(n)t^n\\
G(t)&=&\sum_n g(n)t^n
\end{array}
$$
o\`u $\forall n$, $f(n),g(n) \in \mathbf{Z}$.
Alors la s\'erie produit $FG$ (i.e. le produit de $F$ et $G$ dans l'anneau $\mathbf{Z}[[t]]$) est exactement :
$$FG(t)=\sum_n f\star g (n) t^n.$$
\end{rem}

\begin{rem}
L'op\'erateur $\star$ est associatif et commutatif et l'\'el\'ement unit\'e est la fonction qui vaut $1$ en $0$ et $0$ ailleurs.
\end{rem}

Le th\'eor\`eme suivant permet de calculer les fonctions $q$ \`a partir des fonctions $p$ : 
\begin{thm}
Soient $r\geqslant 1$ et $(m_1,A_1),\ldots,(m_r,A_r)$ des entiers. Alors~:
$$
q(-,(m_1,A_1),\ldots,(m_r,A_r))=\pi_{(m_1,A_1)}\star \ldots    \star\pi_{(m_r,A_r)}
$$
c'est-\`a-dire la distribution associ\'ee \`a $q(-,(m_1,A_1),\ldots,(m_r,A_r))$ est le produit de convolution des distributions associ\'ees \`a $p(-,m_i,A_i)$.
\end{thm}
\begin{proof}
Fixons un $n$.
Pour $r=1$, c'est clair, puisque $q(n,(m_1,A_1))=p(n,m_1,A_1)$.
En effet, $\pi_{(m_1,A_1)}(n_1)=p(n_1,m_1,A_1)$ est exactement le nombre de solutions de : 
$$
\left\lbrace
\begin{array}{rcl}
A_1 \geqslant y_1\geqslant\ldots\geqslant y_m &\geqslant&0,\\
y_1+\ldots+y_m&=&n_1.
\end{array}
\right.
$$

Supposons donc $r \geqslant 2$.
$q(n,(m_1,A_1),\ldots,(m_r,A_r))$ compte le nombre de solutions de : 
$$
\left\lbrace
\begin{array}{l}
   \left\lbrace
   \begin{array}{rcl}
   A_1\geqslant y_{1,1} \geqslant \ldots \geqslant y_{1,m_1} &\geqslant& 0\\
   y_{1,1}+\ldots+y_{1,m_1}&=&n_1
   \end{array}
   \right.
  \\
  \vdots \\
   \left\lbrace
    \begin{array}{rcl}
   A_r\geqslant y_{r,1} \geqslant \ldots \geqslant y_{r,m_r} &\geqslant& 0\\
   y_{r,1}+\ldots+y_{r,m_r}&=&n_r
   \end{array}
   \right.
  \\
  n_1,\ldots,n_r\geqslant 0\\
  n_1+\ldots + n_r=n
\end{array}
\right.
$$
autrement dit, cet entier compte les solutions de~:
$$
\left\lbrace
\begin{array}{l}
\left\lbrace
\begin{array}{l}
   \left\lbrace
   \begin{array}{rcl}
   A_1\geqslant y_{1,1} \geqslant \ldots \geqslant y_{1,m_1} &\geqslant& 0\\
   y_{1,1}+\ldots+y_{1,m_1}&=&n_1
   \end{array}
   \right.
  \\
  \vdots \\
   \left\lbrace
    \begin{array}{rcl}
   A_r\geqslant y_{r-1,1} \geqslant \ldots \geqslant y_{r,m_{r-1}} &\geqslant& 0\\
   y_{r-1,1}+\ldots+y_{r-1,m_{r-1}}&=&n_r
   \end{array}
   \right.
  \\
  n_1,\ldots,n_{r-1}\geqslant 0\\
  n_1+\ldots + n_{r-1}=m
\end{array}
\right.
\\
\left\lbrace
\begin{array}{l}
   \left\lbrace
   \begin{array}{rcl}
   A_r\geqslant y_{1,r} \geqslant \ldots \geqslant y_{1,m_r} &\geqslant& 0\\
   y_{1,r}+\ldots+y_{1,m_r}&=&n_r
   \end{array}
   \right.
   \\
   n_r \geqslant 0\\
\end{array}
\right.

\\
m+n_{r-1}=n
\end{array}
\right.
$$
Donc 
$$
\begin{array}{l}
q(n,(m_1,A_1),\ldots,(m_r,A_r))=\\
\sum_{n_r=1}^n \sum_{m=1}^{n} \pi_{(m_r,A_r)}(n_r) 
\times q(m,(m_1,A_1),\ldots,(m_{r-1},A_{r-1}))\times \chi_n(m+n_r)
\end{array}
$$
o\`u  $\chi_n$ est la fonction caract\'eristique de $n$ qui vaut toujours 0 sauf en $n$ o\`u elle vaut $1$.
Le r\'esultat est alors acquis par r\'ecurrence sur $r$.
\end{proof}

\begin{rem}
Dans la suite, on \'ecrira $y=(A\geqslant y_1\geqslant \ldots \geqslant y_m\geqslant 0)$ pour 
d\'esigner une solution du probl\`eme suivant :
$$
\left\lbrace
\begin{array}{rcl}
A \geqslant y_1\geqslant\ldots\geqslant y_m &\geqslant& 0\\
y_1+\ldots+y_m&=&n\\
\end{array}
\right.
$$
Et on \'ecrira $|y|$ pour $y_1+\ldots + y_m$. 
Donc, $y=(A\geqslant y_1 \geqslant \ldots \geqslant y_m\geqslant 0)$ est une $(m,A)$-partition de $|y|$.
\end{rem}

\Subsection{Fibration en grassmanniennes.}

\begin{defi}
Soit $\lambda=(\lambda_1\geqslant\ldots\geqslant\lambda_d)$ une partition. Soient $c_i$ pour $i\geqslant 1$  des ind\'etermin\'ees.
On d\'efinit le polyn\^ome de Schur correspondant :
$$
\Delta_\lambda(c)=\Delta_{\lambda_1,\ldots,\lambda_d}(c)=\det(c_{\lambda_j+j-i}) \in \mathbf{Z}[c_i,i\in \mathbf{N}^*]
$$
en posant $c_n=0$ si $n<0$, et $c_0=1$.
\end{defi}

\begin{nota}
Soit $\mathcal{E}$ un fibr\'e vectoriel de rang $n$ sur une vari\'et\'e $X$. Soit $d$ un entier plus petit que $n$, posons $Y=\Drapeau{d}{\mathcal{E}}$, le fibr\'e en $d$-grassmanniennes associ\'e \`a $\mathcal{E}$. Notons $f$ la projection $Y\to X$.
On dispose d'une suite exacte universelle : 
$$
\xymatrix{ 
0\ar[r]&\mathcal{S}\ar[r]&f^*(\mathcal{E})\ar[r]&\mathcal{Q}\ar[r]&0 }
$$
o\`u $\mathcal{S}$ est le sous-fibr\'e tautologique de rang $d$, et $\mathcal{Q}$ le fibr\'e quotient tautologique de rang $n-d$.
On pose pour $i\geqslant 0$~: $$
c_i=c_i(\mathcal{Q}-f^*(\mathcal{E}))=c_i(\mathcal{S})
$$
o\`u la fonction $c_i$ est la $i$-\`eme classe de Chern (\cite[Section 3.2]{fulton}  pour leur d\'efinition et propri\'et\'es).
Pour une partition $\lambda=(\lambda_1\geqslant\ldots,\geqslant\lambda_d)$, on notera si aucune ambigu\"it\'e n'est \`a craindre $\Delta_\lambda$ le polyn\^ome de Schur (d\'efinition pr\'ec\'edente) associ\'e \`a ces valeurs de $c_i$.

\end{nota}

D'apr\`es \cite[Proposition 14.6.5]{fulton}, on a le r\'esultat suivant (dit th\'eor\`eme de la base) qui est une g\'en\'eralisation de la formule de Giambella (\cite[Proposition 14.6.4]{fulton}) :
\begin{prop}\label{base1}
Si $X$ est une vari\'et\'e sur $k$ de dimension $N$, et $\mathcal{E}$ un fibr\'e de rang $n$, alors pour $d < n$ et $k \leqslant d(n-d)$,
il existe un isomorphisme canonique
$$
\begin{array}{ccc}
\bigoplus_\lambda \chow^{k-|\lambda|}(X)&\to&\chow^k(\Flag{d}{\mathcal{E}})\\
\alpha_\lambda&\mapsto&\Delta_\lambda . f^*(\alpha_\lambda)
\end{array}
$$
o\`u la somme est prise sur toutes les partitions $\lambda=(\lambda_1 \geqslant \ldots \geqslant \lambda_d)$ telles que $n-d\geqslant \lambda_1 \geqslant \ldots \geqslant \lambda_d \geqslant 0$, et $\Delta_\lambda$ a \'et\'e d\'efinie pr\'ec\'edemment.
\end{prop}

\begin{cor}
Avec les notations de la proposition pr\'ec\'edente, on dispose d'un isomorphisme~:
$$\chow^k(\Drapeau{d}{\mathcal{E}})
\to \bigoplus_{i=0}^{k}\chow^{k-i}(X)^{p(i,d,n-d)}$$
o\`u la fonction $p$ a \'et\'e d\'efinie au paragraphe pr\'ec\'edent.
\end{cor}
\begin{proof}
Il suffit de remarquer que le nombre de partitions $\lambda=(\lambda_1 \geqslant \ldots \geqslant \lambda_d)$ avec $n-d\geqslant \lambda_1 \geqslant \ldots \geqslant \lambda_d \geqslant 0$ avec $|\lambda|=i$ pour un $i$ donn\'e est exactement le nombre de $(d,n-d)$-partitions de $i$.
\end{proof}

\begin{lem}\label{lem_grass}
Soient $S$ une $k$-vari\'et\'e et $\mathcal{E}$ un fibr\'e vectoriel de rang $n$ sur $S$. Soient $r\geqslant 1$ et $1 \leqslant i_1\leqslant\ldots i_{r+1}\leqslant n$ des entiers. Posons $Y=\Drapeau{i_1,\ldots,i_{r+1}}{\mathcal{E}}$ et $X=\Drapeau{i_1,\ldots,i_r}{\mathcal{E}}$.
Alors il existe un fibr\'e vectoriel $\mathcal{W}$ sur $X$ de rang $n-i_r$ tel que l'on ait le diagramme commutatif suivant : 
$$
\xymatrix{
\Drapeau{i_1,\ldots,i_{r+1}}{\mathcal{E}}\ar[rd]\ar@/_/[rdd]\ar[rr]^-\sim&&\Drapeau{i_{r+1}-i_r}{\mathcal{W}}\ar[dl]\ar@/^/[ddl]\\
&\Drapeau{i_1,\ldots,i_r}{\mathcal{E}}\ar[d]&\\
&S&}
$$
\end{lem}
\begin{proof}
On dispose sur $\Drapeau{i_r}{\mathcal{E}}$ du fibr\'e tautologique $\mathcal{V}_r$  tel que 
pour tout $W\in\Drapeau{i_r}{\mathcal{E}}$, on ait $(\mathcal{V}_{r})_{W}=W$ (la construction est la m\^eme
que dans la preuve du th\'eor\`eme \ref{thm_effectif}). Notons $\pi_r$ la projection naturelle
$\Drapeau{i_1,\ldots,i_r}{\mathcal{E}}\to\Drapeau{i_r}{\mathcal{E}}$ et $\pi$ la projection
$\Drapeau{i_1,\ldots,i_r}{\mathcal{E}}\to S$. On a une inclusion $\pi_r^*\mathcal{V}_r \to \pi^*\mathcal{E}$. Alors on construit le fibr\'e vectoriel $\mathcal{W}$ sur $\Drapeau{i_1,\ldots,i_r}{\mathcal{E}}$ comme \'etant le quotient, i.e. tel que la suite suivante soit exacte : 
$$
\xymatrix{
0\ar[r]&\pi_r^* \mathcal{V}_r\ar[r]&\pi^*\mathcal{E}\ar[r]&\mathcal{W}\ar[r]&0.}
$$

Ainsi, pour tout $(V_1,\ldots,V_r) \in X,\, \mathcal{W}_{(V_1,\ldots,V_r)}=\mathcal{E}_{\pi(V_1,\ldots,V_r)}/V_r$. $\mathcal{W}$ est donc un fibr\'e de rang $n-i_r$.

L'application $(V_1,\ldots,V_{r+1})\mapsto V_{r+1}/V_r$ fournit donc pour toute $k$-alg\`ebre $R$ une application : 
$\Drapeau{i_1,\ldots,i_{r+1}}{\mathcal{E}}(R) \to \Drapeau{i_{r+1}-i_r}{\mathcal{W}}(R)$.
Cette application correspond \`a un isomorphisme~:
$$
\xymatrix{
\Drapeau{i_1,\ldots,i_{r+1}}{\mathcal{E}} \ar[r]^-\sim& \Drapeau{i_{r+1}-i_r}{\mathcal{W}}. }\qed
$$\noqed
\end{proof}

\Subsection{Fibration en vari\'et\'es de drapeaux}

En it\'erant le lemme \ref{lem_grass}, on obtient la proposition suivante~:
\begin{prop}\label{base}
Si $S$ est une vari\'et\'e sur $k$, $\mathcal{E}$ un fibr\'e de rang $n$ sur $S$, $r\geqslant 1$ un entier, et $1 \leqslant i_1<\ldots<i_r \leqslant n$ une suite d'entiers, alors pour tout $k$ , on a un isomorphisme~:
$$
\chow^k(\Drapeau{i_1,i_2,\ldots,i_r}{\mathcal{E}}) \to
\bigoplus_{\lambda_1,\ldots,\lambda_r} \chow^{k-|\lambda_1|-\ldots-|\lambda_r|}(S)
$$
o\`u la somme est prise sur toutes les partitions $\lambda_s=(n-i_s\geqslant \lambda_{s,1} \geqslant \ldots \geqslant \lambda_{s,i_s-i_{s-1}} \geqslant 0)$ pour $s\in \lbrace 1,\ldots,r\rbrace$, en d\'efinissant $i_0=0$.
\end{prop}
\begin{proof}
Posons $X_0=S$, et pour $s\in \lbrace 1,\ldots,n\rbrace$, $X_s=\Drapeau{i_1,\ldots,i_s}{\mathcal{E}}$. \\
D'une part, $X_1 \to X_0$ \'etant par construction m\^eme une fibration en grassmanniennes, on a, d'apr\`es la proposition 
\ref{base1}, un isomorphisme~:
$$
\chow^k(X_{1})\to \bigoplus_{\lambda_1}\chow^{k-|\lambda_1|}(X_0)
$$
o\`u la somme est prise sur toutes les partitions $\lambda_1=(\lambda_{1,1} \geqslant \ldots \geqslant \lambda_{1,i_1})$ telles que $n-i_1\geqslant \lambda_{1,1} \geqslant \ldots \geqslant \lambda_{1,i_1} \geqslant 0$.\\
D'autre part, d'apr\`es le lemme \ref{lem_grass}, pour tout $s\in \lbrace 2\ldots r \rbrace$, la projection naturelle $X_{s} \to X_{s-1}$ est une fibration en grassmanniennes, i.e. : il existe un fibr\'e vectoriel $\mathcal{W}_s$ sur $X_{s-1}$ de rang $n-i_{s-1}$ tel que le diagramme suivant soit commutatif :
$$
\xymatrix{
X_{s}\ar[rd]\ar[rr]^-\sim&&\Drapeau{i_{s}-i_{s-1}}{\mathcal{W}_s}\ar[dl]\\
&X_{s-1}&}
$$
et ainsi d'apr\`es la proposition \ref{base1}, on a un isomorphisme~: 
$$
\chow^k(X_{s})\to \bigoplus_{\lambda_s}\chow^{k-|\lambda^s|}(X_{s-1})
$$
o\`u la somme est prise sur toutes les partitions $\lambda_s=(\lambda_{s,1}\geqslant \ldots \geqslant \lambda_{s,i_s-i_{s-1}})$ telles que $n-i_s\geqslant \lambda_{s,1} \geqslant \ldots \geqslant \lambda_{s,i_s-i_{s-1}} \geqslant 0$.

Par r\'ecurrence, il suffit alors de consid\'erer le diagramme commutatif suivant :
$$
\xymatrix{
X_r\ar[d]\ar[r]^-\sim&\Drapeau{i_{r}-i_{r-1}}{\mathcal{W}_{r}}\ar[dl]\\
X_{r-1}\ar[d]\ar[r]^-\sim&\Drapeau{i_{r-1}-i_{r-2}}{\mathcal{W}_{r-1}}\ar[dl]\\
X_{r-2}\ar[d]&\\
\vdots\ar[d]&\\
X_1\ar[d]\ar[r]^-\sim&\Drapeau{i_1}{\mathcal{W}_1}\ar[dl]\\
X_0&}
$$
En particulier $\mathcal{W}_1=\mathcal{E}$.
Donc on obtient un isomorphisme~:
$$
\chow^k(\Drapeau{i_1,i_2,\ldots,i_r}{\mathcal{E}}) \to
\bigoplus_{\lambda_1,\ldots,\lambda_r} \chow^{k-|\lambda_1|-\ldots-|\lambda_r|}(S)
$$
o\`u la somme est prise sur toutes les partitions $\lambda_s=(\lambda_{s,1}\geqslant \ldots \geqslant \lambda_{s,i_s-i_{s-1}})$ telles que $n-i_s\geqslant \lambda_{s,1} \geqslant \ldots \geqslant \lambda_{s,i_s-i_{s-1}} \geqslant 0$ pour $s=1,\ldots,r$, en d\'efinissant $i_0=0$.
\end{proof}

\begin{cor}\label{cor_coeff}
Conservons les notations de la proposition pr\'ec\'edente, on a un isomorphisme~: 
$$
\chow^k(\Drapeau{i_1,i_2,\ldots,i_r}{\mathcal{E}}) \to
\bigoplus_{i=0}^{k}(\chow^{k-i}(S))^{n_i}
$$
o\`u $n_i=q(i,(i_1,n-i_1),(i_2-i_1,n-i_2),\ldots,(i_r-i_{r-1},n-i_r))$ (la fonction $q$ est d\'efinie au paragraphe pr\'ec\'edent).
\end{cor}
\begin{proof}
Il suffit de remarquer que pour un $i$ donn\'e, le nombre d'occurrences de $i$ dans la somme
 $|\lambda_1|+\ldots+|\lambda_r|$  pour toutes les partitions 
$\lambda_s=(\lambda_{s,1}\geqslant \ldots \lambda_{s,i_s-i_{s-1}})$ telles que
 $n-i_s\geqslant \lambda_{s,1} \geqslant \ldots \geqslant \lambda_{s,i_s-i_{s-1}} \geqslant 0)$ 
pour $s\in \lbrace 1,\ldots,r\rbrace$ est exactement le nombre de
$((i_1,n-i_1),(i_2-i_1,n-i_2),\ldots,(i_r-i_{r-1},n-i_r))$-partitions de $i$. Donc il y a exactement $q(i,(i_1,n-i_1),(i_2-i_1,n-i_2),\ldots,(i_r-i_{r-1},n-i_r))$ occurrences de $i$ dans la somme.
\end{proof}

\Subsection{Fibration en produit de vari\'et\'es de drapeaux.}
Ici encore, il s'agit d'une cons\'equence du lemme \ref{lem_grass}.

\begin{prop}\label{base_prod}
Soient $S$ une vari\'et\'e sur $k$, $\mathcal{E}_{-}$ et  $\mathcal{E}_{+}$   deux fibr\'es vectoriels sur $S$ de rang $n_{-}$ et $n_{+}$ respectivement, $1\leqslant s \leqslant r \leqslant n$ deux  entiers, et $1 \leqslant i_1<\ldots<i_s \leqslant n_-$  et $1 \leqslant i_{s+1}<\ldots<i_r \leqslant n_+$ deux suites d'entiers.
Posons $X=\Drapeau{i_1,\ldots,i_s}{\mathcal{E}_{-}}$ 
et $Y=\Drapeau{i_{s+1},\ldots,i_r}{\mathcal{E}_{+}}$.
Alors pour tout $k$, on a un isomorphisme~:
$$
\chow^k(X\times_S Y) \to 
\bigoplus_{\lambda_1,\ldots,\lambda_r} \chow^{k-|\lambda_1|-\ldots-|\lambda_r|}(S)
$$
o\`u la somme est prise sur toutes les partitions : 
$$
\begin{array}{rcl}
\lambda_1&=&(n_{-}-i_1\geqslant \lambda_{1,1}\geqslant\ldots \geqslant \lambda_{1,i_1}\geqslant 0)\\
\lambda_2&=&(n_{-}-i_2\geqslant \lambda_{2,1}\geqslant\ldots  \geqslant \lambda_{2,i_2-i_1}\geqslant 0)\\
\vdots&\vdots&\vdots\\
\lambda_s&=&(n_{-}-i_{s}\geqslant \lambda_{s,1}\geqslant\ldots  \geqslant \lambda_{s,i_{s}-i_{s-1}}\geqslant 0)\\
\lambda_{s+1}&=&(n_{+}-i_{s+1}) \geqslant \lambda_{s+1,1}\geqslant \ldots  \geqslant \lambda_{s+1,i_{s+1}} \geqslant 0)\\
\lambda_{s+2}&=&(n_{+}-i_{s+2} \geqslant \lambda_{s+2,1}\geqslant \ldots  \geqslant \lambda_{s+2,i_{s+2}-i_{s+1}} \geqslant 0)\\
\vdots&\vdots&\vdots\\
\lambda_{r}&=&(n_{+}-i_{r} \geqslant \lambda_{r,1}\geqslant \ldots  \geqslant \lambda_{r,i_{r}-i_{r-1}} \geqslant 0).
\end{array}
$$
\end{prop}

\begin{proof}
Consid\'erons le carr\'e cart\'esien suivant : 
$$
\xymatrix{
  X\times_S Y
  \ar[r]^-{\pi_X}
  \ar[d]^-{\pi_Y}
  \cartesien
  &
  X
  \ar[d]^-{f_X}
  \\
  Y
  \ar[r]^-{f_Y}
  &
  S}
$$
Alors il suffit d'appliquer la proposition \ref{base} aux deux fibrations en vari\'et\'es de drapeaux : 
$X\to S$ et $X\times_S Y \to X$.
\end{proof}

\begin{cor}\label{cor_coeff_prod}
Conservons les notations de la proposition pr\'ec\'edente. On a un isomorphisme~: 
$$
\chow^k(X \times_S Y)\to
\oplus_{i=0}^{k}(\chow^{k-i}S)^{n_i}
$$
o\`u 
$n_i=
q(i,(i_1,n_{-}-i_1),(i_2-i_1,n_{-}-i_2),\ldots,(i_s-i_{s-1},n_{-}-i_{s}),
(i_{s+1},n_{+}-i_{s+1}),(i_{s+2}-i_{s+1},n_{+}-i_{s+2}),\ldots,(i_r-i_{r-1},n_{+}-i_r)).$
(la fonction $q$ a \'et\'e d\'efinie au paragraphe pr\'ec\'edent).
\end{cor}
\begin{proof}
La preuve est analogue \`a celle du corollaire \ref{cor_coeff}.
\end{proof}

\section{Applications}

\Subsection{Groupes de Chow des vari\'et\'es de drapeaux tordues}

Gr\^ace aux r\'esultats sur le calcul des groupes de Chow des vari\'et\'es de drapeaux du paragraphe pr\'ec\'edent, au th\'eor\`eme \ref{thm_effectif} et aux corollaires \ref{cor_thm_general} et \ref{cor_ind_premier}, on est en mesure de calculer de fa\c{c}on explicite les groupes de Chow des vari\'et\'es tordues sous quelques hypoth\`eses.
\Subsubsection{Cas $i_1=1$}
Si l'on est sous les hypoth\`eses du  th\'eor\`eme \ref{thm_effectif} : 
\begin{prop}\label{prop_coeff_thm_effectif}
Soient $A$ une alg\`ebre d'Azumaya sur $k$ de degr\'e $n$ et $1\leqslant i_1<\ldots<i_r \leqslant n$ des entiers. On suppose que $i_1=1$. 
Alors on a un isomorphisme~:
$$\chow^k(\Flag{1,i_2,\ldots,i_r}{A})\to
\bigoplus_{\lambda_r,\ldots,\lambda_2} \chow^{k-|\lambda_r|-\ldots-|\lambda_2|}(\SB{}{A})$$
o\`u la somme est prise sur toutes les partitions
 $\lambda_s=(\lambda_{s,1}\geqslant \ldots \geqslant \lambda_{s,i_{r-s+3}-i_{r-s+2}})$ telles que $i_{r-s+2}-1 \geqslant \lambda_{s,1} \geqslant \ldots \geqslant \lambda_{s,i_{r-s+3}-i_{r-s+2}} \geqslant 0$ pour $s\in \lbrace 2,\ldots,r\rbrace$ en posant $i_{r+1}=n$,
et un isomorphisme~:
$$
\chow^k(\Flag{1,i_2,\ldots,i_r}{A})\to
\bigoplus_{i=0}^{k}\chow^{k-i}(\SB{}{A})^{n_{i}}
$$
avec $n_{i}=q(i,(i_{r+1}-i_r,i_{r}-1),(i_{r}-i_{r-1},i_{r-1}-1),\ldots,(i_3-i_2,i_2 -1 ))$.      
\end{prop}

\begin{proof}
C'est l'application de la proposition \ref{base} et du corollaire \ref{cor_coeff} \`a un fibr\'e de rang $n-1$ sachant le r\'esultat du th\'eor\`eme \ref{thm_effectif}.
\end{proof}

\Subsubsection{Sous les hypoth\`eses du corollaire \ref{cor_thm_general}}
%Si l'on est sous les hypoth\`eses du corollaire \ref{cor_thm_general} :
\begin{prop}\label{prop_coeff_thm_general}
Soient $A$ une alg\`ebre d'Azumaya sur $k$ de degr\'e $n$ et $1\leqslant i_1<\ldots<i_r \leqslant n$ des entiers. Conservons les hypoth\`eses du corollaire \ref{cor_thm_general}.
Alors on a un isomorphisme~: 
$$
\chow^k(\Flag{i_1,i_2,\ldots,i_r}{A})\to
\bigoplus_{\lambda_1,\ldots,\widehat{\lambda_s},\ldots,\lambda_r}
\chow^{k-|\lambda_1|-\ldots-\widehat{|\lambda_s|}-\ldots-|\lambda_r|}(\SB{i_s}{A})
$$
o\`u la somme est prise sur toutes les partitions
$$
\begin{array}{rcl}
\lambda_1&=&(i_s-i_1\geqslant \lambda_{1,1}\geqslant\ldots \geqslant \lambda_{1,i_1}\geqslant 0)\\
\lambda_2&=&(i_s-i_2\geqslant \lambda_{2,1}\geqslant\ldots  \geqslant \lambda_{2,i_2-i_1}\geqslant 0)\\
\vdots&\vdots&\vdots\\
\lambda_{s-1}&=&(i_s-i_{s-1}\geqslant \lambda_{s-1,1}\geqslant\ldots  \geqslant \lambda_{s-1,i_{s-1}-i_{s-2}}\geqslant 0)\\
\lambda_{s+1}&=&(n-i_{s+1}\geqslant \lambda_{s+1,1}\geqslant \ldots  \geqslant \lambda_{s+1,i_{s+1}-i_s} \geqslant 0)\\
\lambda_{s+2}&=&(n-i_{s+2} \geqslant \lambda_{s+2,1}\geqslant \ldots  \geqslant \lambda_{s+2,i_{s+2}-i_{s+1}} \geqslant 0)\\
\vdots&\vdots&\vdots\\
\lambda_{r}&=&(n-i_{r} \geqslant \lambda_{r,1}\geqslant \ldots  \geqslant \lambda_{r,i_{r}-i_{r-1}} \geqslant 0).
\end{array}
$$
De plus, on a un isomorphisme~: 
$$
\chow^k(\Flag{i_1,i_2,\ldots,i_r}{A})\to
\bigoplus_{i=0}^{k}\chow^{k-i}(\SB{i_s}{A})^{n_i}
$$
o\`u $n_i=
q(i,(i_1,i_s-i_1),(i_2-i_1,i_s-i_2),\ldots,(i_{s-1}-i_{s-2},i_{s}-i_{s-1}),
(i_{s+1}-i_s,n-i_{s+1}),(i_{s+2}-i_{s+1},n-i_{s+2}),\ldots,(i_r-i_{r-1},n-i_r)).$
\end{prop}

\begin{proof}
C'est l'application de la proposition \ref{base_prod} et du corollaire \ref{cor_coeff_prod} sachant les r\'esultats du th\'eor\`eme \ref{thm_general} et du corollaire \ref{cor_thm_general}.
\end{proof}

\begin{rem}\label{rem_premier}
Si l'on est sous les hypoth\`eses du corollaire \ref{cor_ind_premier}, alors, avec les m\^emes notations, on dispose d'un $s \in \lbrace 1,\ldots,r \rbrace$ tel que le corollaire \ref{cor_thm_general} puisse s'appliquer. Et par suite, la proposition pr\'ec\'edente est encore valable avec ce $s$.
\end{rem}

%\begin{rem}
%Soit $A$ une alg\`ebre d'Azumaya sur $k$ de degr\'e $n$. Soient $r>1$ et $1\geqslant i_1<i_2<\ldots<i_r \leqslant n%$ des entiers. Si  $\Flag{i_1,\ldots,i_r}{A}$ v\'erifie les hypoth\'eses du th\'eor\`eme \ref{thm_effectif} ou des %%%%corollaires \ref{cor_thm_general} ou \ref{cor_ind_premier}, on dira que $(A,i_1,\ldots,i_r)$ v\'erifie l'hypoth\'ese %%$(\mathcal{D})$. Et on posera alors $S=SB{i_s}{A}$ avec $i_s=1$ si l'on est dans les hypoth\`eses du th\'eor\`me \re%f{thm_effectif}.
%\end{rem}

Remarquons que l'on a $n_0=1$. Ainsi, dans les deux propositions  pr\'ec\'edentes (\ref{prop_coeff_thm_effectif} et \ref{prop_coeff_thm_general}) , pour un $k$ donn\'e,  il existe des entiers $n_i$ tels que l'on ait un isomorphisme~: 
$$
\chow^k(X)\to\chow^k(S)\oplus(\oplus_{i\geqslant 1} \chow^{k-i}(S)^{n_i})
$$
 o\`u $X$ est la vari\'et\'e de drapeaux consid\'er\'ee et $S$ la base sur laquelle on projette ($\SB{}{A}$ ou $\SB{i_s}{A}$). Et, par suite, on a : 

\begin{cor}
Pla\c{c}ons-nous, soit dans les hypoth\`eses de la proposition \ref{prop_coeff_thm_effectif} soit dans celles de la proposition  \ref{prop_coeff_thm_general}. Notons $X$  la vari\'et\'e de drapeaux consid\'er\'ee et $S$ la base (vari\'et\'e de Severi-Brauer) sur laquelle on projette.
Soit $k$ un entier. Si les groupes de Chow de $S$ n'ont pas de torsion en codimension plus petite que $k$, alors il en va de m\^eme pour ceux de $X$, et la r\'eciproque est vraie.
\end{cor}

\Subsubsection{Calcul \`a partir d'une vari\'et\'e de Severi-Brauer}
Fixons une alg\`ebre d'Azumaya $A$ sur $k$ de degr\'e $n$. Soient $1\leqslant i_1 < \ldots < i_r \leqslant n$ des entiers. Supposons que $(\ind A,i_1,\ldots,i_r)=1$. Posons $X=\Flag{i_1,\ldots,i_r}{A}$ et $S=\SB{}{A}$. On souhaite calculer les groupes de Chow de $X$ \`a partir de ceux de $S$. Si $i_1=1$, alors la proposition \ref{prop_coeff_thm_effectif} permet de voir que les groupes de Chow de $X$ se calculent explicitement \`a partir de ceux de $S$. Supposons donc d\'esormais que $1<i_1$. Posons $Y=\Flag{1,i_1,\ldots,i_r}{A}$. 
Consid\'erons le diagramme suivant : 
$$
\xymatrix{&Y\ar[dl]_-{\pi}\ar[dr]^-{f}&\\
S&&X}
$$
o\`u la fl\`eche $f$ est l'application $(I,I_1,\ldots,I_r) \mapsto (I_1,\ldots,I_r)$.
\begin{lem}
l'application $f\,:\,Y\to X$ est une fibration en espaces projectifs.
\end{lem}
\begin{proof}
L'application $f$ est une fibration en vari\'et\'es de Severi-Brauer. 
Puisque $(\ind A,i_1,\ldots,i_r)=1$, alors l'image de la classe de $A$ dans $\brauer(k(X))$ est triviale (cf \cite[(5.11)]{merkurjev_index_1}). En outre, la classe de $A$ et celle de $\mathcal{B}_-$ ont m\^eme image dans $\brauer(k(X))$. On en d\'eduit alors que la classe de $\mathcal{B}_-$ est triviale.
%Consid\'erons le diagramme suivant, qui est cart\'esien : 
%$$
%\xymatrix{
%Y\ar[d]^-{f}\ar[r]&\Flag{1,i_1}{A}\ar[d]\\
%X\ar[r]_-{g}&\Flag{i_1}{A} }
%$$
%Or $\Flag{1,i_1}{A} \to \Flag{i_1}{A}$ est d'apr\`es le th\'eor\`eme %\ref{thm_general} une fibration en vari\'et\'e de Severi-Brauer. Comme l'alg\`ebre %$\mathcal{B}_{-}$ intervenant dans la fibration est \'egale \`a la classe de $A$ %sur chaque corps r\'esiduel, alors puisque $(\ind A,i_1,\ldots,i_r)=1$, %$g^*(\mathcal{B}_{-})$ est d\'eploy\'ee sur $X$. 
\end{proof}

On a donc, pour un entier $k$, les deux isomorphismes suivants~: 
$$
\chow^k(Y)\to\oplus_{i=1}^{n-1}\chow^{k-i}(X) \qquad\text{et}\qquad
\chow^k(Y)\to\oplus_{i}\chow^{k-i}(S)^{n_{i}}
$$
avec $n_{i}=q(i,(n-i_r,i_r-1),(i_r-i_{r-1},i_{r-1}-1),\ldots,(i_2-i_1,i_1-1))$. La seconde \'equation r\'esulte simplement de  la proposition 
 \ref{prop_coeff_thm_effectif}. Quant \`a la premi\`ere \'equation, il s'agit de la formule classique du calcul des groupes de Chow d'un fibr\'e projectif en fonction des groupes de Chow de sa base (cas particulier de la proposition \ref{base1}).

On a donc la proposition suivante :
\begin{prop}
Pour tout $k$ on a une suite exacte :
$$
\xymatrix{
0\ar[r]&\oplus_{i}\chow^{k-1-i}(S)^{n_{i}}\ar[r]&\oplus_{i}\chow^{k-i}(S)^{n_{i}}\ar[r]&\chow^{k}(X)\ar[r]&0.}
$$
\end{prop}
\begin{proof}
On a en effet une suite exacte pour tout $k$ :
$$
\xymatrix{
0\ar[r]&\chow^{k-1}(Y)\ar[r]&\chow^k(Y)\ar[r]&\chow^{k}(X)\ar[r]&0}.
$$
\end{proof}
\begin{rem}
Ainsi $\chow^k(X)$ s'identifie \`a un quotient du groupe $\oplus_{i}\chow^{k-i}(S)^{n_{i}}$
par $\oplus_{i}\chow^{k-1-i}(S)^{n_{i}}$, mais l'injection 
$$\oplus_{i}\chow^{k-1-i}(S)^{n_{i}}\to \oplus_{i}\chow^{k-i}(S)^{n_{i}}$$ n'est pas simple et fait intervenir les polyn\^omes de Schur li\'es \`a la fibration $Y \to X$.
\end{rem}
On peut en d\'eduire ais\'ement le corollaire suivant :
\begin{cor}
Conservons les m\^emes notations.
Soit $k$ un entier. Si les groupes de Chow de $S$ n'ont pas de torsion en codimension plus petites que $k$, alors il en va de m\^eme pour ceux de $X$, et la r\'eciproque est vraie.
\end{cor}

\Subsection{Exemple : $\chow^2$ sans torsion}

\begin{rem}
Supposons que $X\to S$ soit une projection de $k$-vari\'et\'es telle qu'il existe des entiers $n_{k,i}$, de sorte que pour tout $k$, on ait un isomorphisme~: 
$$
\chow^k(X)\to\oplus_{i=0}^{k}(\chow^{k-i}(S))^{n_{k,i}}.
$$
Alors, si les groupes de Chow de $S$ en codimension $\leqslant k_0$ n'ont pas de torsion, alors les groupes de Chow de $X$ en codimension $\leqslant k_0$ n'ont pas de torsion.
\end{rem}

\begin{prop}
Si $X$ est une vari\'et\'e projective homog\`ene sous un groupe lin\'eaire, alors $\chowtors{0}{X} =0$ et $\chowtors{1}{X}=0$.
\end{prop}
\begin{proof}
Pour le $\CH^0$ le r\'esultat est clair. Pour le $\CH^1$, il suffit de remarquer que $\CH^1 X = \Pic X$. Or d'apr\`es \cite[lemme 5.1]{peyre_galois3}, on a une suite exacte :
$$
\xymatrix{0\ar[r]& \Pic X \ar[r] &\Pic (X_{k_s})^\Gamma }
$$
o\`u $k_s$ est une cl\^oture s\'eparable de $k$ et $\Gamma$ le groupe de Galois absolu de $k$. Ainsi, $\CH^1 X$ est sans torsion.
\end{proof}

D'apr\`es Karpenko (\cite{karpenko_filtration}), on a :
\begin{prop}[Karpenko]
Si $A$ est une alg\`ebre d'Azumaya sur un corps $k$, dont l'indice co\"incide avec son
exposant, alors $$\chowtors{2}{\SB{}{A}}=0.$$
\end{prop}

\begin{prop}
Si $A$ est une alg\`ebre d'Azumaya de degr\'e $n$, dont l'indice co\"incide avec son
exposant et $1= i_1<\ldots<i_r \leqslant n$ alors $\chowtors{2}{\Flag{i_1,i_2,\ldots,i_r}{A}}=0$.
\end{prop}
\begin{proof}
Sachant que pour $i\leqslant 2$, $\chowtors{i}{SB(A)}=0$, il suffit d'appliquer la remarque du d\'ebut de section.
\end{proof}

\Subsection{Exemple : vari\'et\'es de drapeaux complets et torsion dans les $\chow^i$, $i\geqslant 2$}

Commen\c{c}ons par \'enoncer un r\'esultat de Karpenko : 
\begin{prop}[Karpenko]
Soit $e$ un entier non divisible par un carr\'e. Alors il existe une alg\`ebre d'Azumaya $A$ d'exposant $e$ dont $\chowtors{2}{\SB{}{A}}$ est d'ordre $e$.
\end{prop}
\begin{proof}
Corollaire 5.2 du chapitre 1 de \cite{karpenko_hab}.\\
\end{proof}
\begin{rem}
 Dans cette proposition, on peut en fait choisir $A$ d'indice (et donc de degr\'e) aussi grand que l'on veut (proposition 5.1 du chapitre 1 de \cite{karpenko_hab}). La construction d'une telle alg\`ebre est faite dans l'exemple 4.12 du chapitre 1 de \cite{karpenko_hab}.
\end{rem}

Examinons le cas des vari\'et\'es de drapeaux complets.
Soit $A$ une alg\`ebre d'Azumaya sur $k$ de degr\'e $n$. Notons $X=\Flag{1,\ldots,n}{A}$ la vari\'et\'e de drapeaux complets. 
Alors, gr\^ace au th\'eor\`eme \ref{thm_effectif}, on dispose d'un fibr\'e $\mathcal{V}$ de rang $n-1$ sur $\SB{}{A}$ tel que 
$$
\xymatrix{\Flag{1,2,\ldots,n}{A}\ar[rd]\ar[rr]^\sim&&\Drapeau{1,\ldots,n-1}{\mathcal{V}}\ar[ld]\\
&\SB{}{A}&}
$$
et donc par la proposition \ref{prop_coeff_thm_effectif}, on a un isomorphisme~: 

$$\chow^k(X)\to
\oplus_{i=0}^{k}(\chow^{k-i}\SB{}{A})^{n_i}
$$
avec $n_i=q(i,(1,n-1),\ldots,(1,2),(1,1))$. L'entier $n_i$ compte le nombre de solutions de : 
$$
\left\lbrace
\begin{array}{rcl}
\forall j \,\, x_j \leqslant j\\
x_1+\ldots+x_{n-1}&=&i\\
\end{array}
\right.
$$
donc pour tout $i \in \lbrace 0,\ldots, \frac{(n-1)n}{2} \rbrace$, $n_i \geqslant 1$. Ainsi :

\begin{cor}
Soit $e$ un entier non divisible par un carr\'e. Pour tout entier $k\geqslant 2$, il existe une alg\`ebre d'Azumaya $A$ d'exposant $e$ et une vari\'et\'e homog\`ene $X$ sous $\PGL(A)$, telles que $\chow^k(X)$ contienne un sous-groupe cyclique d'ordre $e$.
\end{cor}
\begin{proof}
 Supposons que  $A$ soit une alg\`ebre d'Azumaya d'exposant $e$ telle que $\chowtors{2}{\SB{}{A}}$ soit un groupe cyclique d'ordre $e$. Si on pose $X=\Flag{1,2,\ldots,n}{A}$ o\`u $n$ est le degr\'e de $A$, alors d'apr\`es ce qui pr\'ec\`ede, $\chow^k(X)$ contiendra n\'ecessairement un sous-groupe cyclique d'ordre $e$ d\`es qu'il existera un $i$ satisfaisant :
$$\left \lbrace
\begin{array}{l}
k-i=2\\
i \leqslant \frac{(n-1)n}{2}.
\end{array} \right.
$$
Un tel $i$ existe si et seulement si :
$$ 
k \leqslant \frac{(n-1)n}{2} + 2.
$$
Pour un $k$ fix\'e, il suffit alors de prendre une alg\`ebre d'Azumaya $A$ d'exposant $e$ donn\'ee par la proposition pr\'ec\'edente et de degr\'e $n$ assez grand.
\end{proof}

\nocite{*}
\bibliographystyle{./franckplain}
\bibliography{./art-drap}

\end{document}